\documentclass[a4paper,12pt]{article}
\pdfoutput=1
\usepackage[cp1251]{inputenc}
\usepackage[T2A]{fontenc}
\usepackage[english]{babel}
\usepackage[intlimits]{amsmath}
\usepackage{amsthm}
\usepackage{bm}
\usepackage{amsfonts,amssymb}

\newcommand{\supp}{\mathop{\text{supp}}}

\def\leq{\leqslant}
\def\geq{\geqslant}

\theoremstyle{plain}
\newtheorem{remark}{Remark}
\newtheorem{lemma}{Lemma}
\newtheorem{prop}{Proposition}
\newtheorem{theorem}{Theorem}
\newtheorem{corollary}{Corollary}

\begin{document}
\title{The effect of curvature \\ in fractional Hardy--Sobolev inequality \\ with the Spectral Dirichlet Laplacian}

\author{Nikita Ustinov\footnote{St.Petersburg State University, 7/9 Universitetskaya Emb., St.Petersburg 199034, Russia. E-mail: {ustinns@yandex.ru}.  Supported by RFBR grant 17-01-00678A}}

\date{}

\maketitle

{\footnotesize
\noindent
{\bf Abstract.} We prove the attainability of the best constant in the fractional Hardy--Sobolev inequality with boundary singularity for the Spectral Dirichlet Laplacian. The main assumption is the average concavity of the boundary at the origin.
\medskip


\medskip
\noindent
{\bf 2010 Mathematics Subject Classfication:} 35R11, 49J10, 35B45}

\section{Introduction}
In this paper we discuss the attainability of the best fractional Hardy--Sobolev constant $\mathcal{S}^{Sp}_{s, \sigma}(\Omega)$ in $\mathcal{C}^1$--smooth bounded domain $\Omega \subset \mathbb{R}^n,$ $n \geq 2:$
\begin{equation}
\label{in1}
\mathcal{S}^{Sp}_{s, \sigma}(\Omega) \cdot \| |x|^{\sigma-s} u  \|^2_{L_{2^*_{\sigma}}(\Omega)} \leq \langle (-\Delta)_{Sp}^s u, u \rangle, \quad u \in \widetilde{\mathcal{D}}^s(\Omega),
\end{equation}
where $0 < \sigma < s < 1$ and $2^*_{\sigma} \equiv \frac{2n}{n-2\sigma}.$ Fractional Laplacian in the right-hand side of~(\ref{in1}) is the \textit{Spectral Dirichlet Laplacian}; the space $\widetilde{\mathcal{D}}^s(\Omega)$ is generated by its quadratic form (see Sec. 2).

In case of $0 \not\in \overline{\Omega}$ the embedding $\widetilde{\mathcal{D}}^s(\Omega) \hookrightarrow L_{2^*_{\sigma}}(\Omega, |x|^{(\sigma-s)2^*_{\sigma}})$ is compact and $\mathcal{S}^{Sp}_{s, \sigma}(\Omega)$ is obviously attained. Through this paper we will consider nontrivial case $0 \in \overline{\Omega}.$

In the local case $s = 1$ the inequality (\ref{in1}) coincides with the inequality 
\begin{equation}
\label{in_local}
\mathcal{S}_{\sigma}(\Omega) \cdot \| |x|^{\sigma-1} u  \|^2_{L_{2^*_{\sigma}}(\Omega)} \leq \langle -\Delta u, u \rangle  = \| \nabla u\|^2_{L_2(\Omega)}.
\end{equation}
Attainability of the best constant $\mathcal{S}_{\sigma}(\Omega)$ is well-studied (even for the non-Hilbertian case):
\begin{itemize}
\item If $0 \in \Omega,$ $\sigma \in [0,1]$ and $n \geq 3,$ then $\mathcal{S}_{\sigma}(\Omega)$ does not depend on $\Omega;$ for~$\sigma \in (0,1]$ the constant $\mathcal{S}_{\sigma}(\mathbb{R}^n)$ \textit{is attained} on the family of functions
\begin{equation*}
u_{\varepsilon}(x) := \left( \varepsilon + |x|^{\frac{2\sigma(n-2)}{n-2\sigma}} \right)^{1-\frac{n}{2\sigma}},
\end{equation*}
(\cite{Glaser, Lieb}; in non-Hilbertian case see \cite{Aubin, Talenti} for $\sigma = 1,$ \cite{Gh_Yu} for $\sigma \in (0,1)$), thus $\mathcal{S}_{\sigma}(\Omega)$ \textit{is not attained} if $\widetilde{\mathcal{D}}^1(\Omega) \neq \mathcal{D}^{1}(\mathbb{R}^n);$ if $\sigma = 0,$ then $\mathcal{S}_{\sigma}$ \textit{is not attained} even in $\mathbb{R}^n$ (see \cite[Sec. 7.3]{Hardy}).
\item In case of $0 \in \partial \Omega$ attainability of $\mathcal{S}_{\sigma}(\Omega)$ was proved for cones: if $\sigma \in (0,1),$ $n \geq 2$ and~$\Omega$ is a cone in $\mathbb{R}^n,$ then $\mathcal{S}_{\sigma}(\Omega)$ \textit{is attained} (\cite{Egn}; \cite{Naz_Cone} in non-Hilbertian case).
\item The case of bounded $\Omega$ with $0 \in \partial \Omega$ is much more complex and depends on the behaviour of $\partial \Omega$ at the origin. In \cite{Gh_Kang} it was shown that for $n \geq 4$ $\mathcal{S}_{\sigma}(\Omega)$ \textit{is attained} if all principal curvatures of $\partial \Omega$ are negative at the origin. In \cite{Gh_Rob} this condition was replaced by the negativity of the mean curvature of $\partial \Omega$ at the origin. In \cite{Dem_Naz} these conditions were sufficiently weakened and the \textit{attainability} was proved for all $n \geq 2.$ 
\end{itemize}
For $s \not \in \mathbb{N}$ only a few results were established before. In \cite{Yang} \textit{attainability} of $\mathcal{S}^{Sp}_{s, \sigma}(\mathbb{R}^n)$ was shown for $s \in \left(0, \tfrac{n}{2}\right).$ For $s \in (0,1)$ \textit{attainability} of the best constant in $\mathbb{R}^n_{+}$ was shown for fractional Hardy--Sobolev inequalities with restricted Dirichlet and Neumann fractional Laplacians \cite{MN1, MN5}. These inequalities differ from (\ref{in1}) by the choice of fractional Laplacian in the right-hand side.

\medskip

In this paper we prove the following results for the inequality (\ref{in1}):
\begin{itemize}
\item In case of $0 \in \Omega$ and $\widetilde{\mathcal{D}}^s(\Omega) \neq \mathcal{D}^{s}(\mathbb{R}^n)$ the best constant $\mathcal{S}^{Sp}_{s, \sigma}(\Omega)$ \textit{is not attained}. Moreover, if the domain $\Omega$ is star-shaped around the origin, then the corresponding Euler--Lagrange equation has \textit{only trivial non-negative solution}.
\item The best constant $\mathcal{S}^{Sp}_{s, \sigma}(\mathbb{R}_+^n)$ \textit{is attained}.
\item In case of $0 \in \partial \Omega$ in bounded $\Omega$ the best constant $\mathcal{S}^{Sp}_{s, \sigma}(\Omega)$ \textit{is attained} under the assumptions on $\partial \Omega$ at the origin, analogous to the conditions from \cite{Dem_Naz}.
\end{itemize}
The short announcement of these results was given in \cite{Ustinov}.
\medskip

The paper consists of nine sections. In Sec. 2 we give basic definitions and recall some properties of the \textit{Spectral Dirichlet Laplacian} (including the Stinga--Torrea extension).
In Sec. 3 we prove \textit{unattainability} of $\mathcal{S}^{Sp}_{s, \sigma}(\Omega)$ in case of $0 \in \Omega$ together with \textit{non-existence of positive solutions} for the Euler--Lagrange equation in star-shaped $\Omega.$ In Sec. 4 we derive estimates for Green functions of some auxiliary problems. In Sec. 5 we prove \textit{attainability} of the best constant $\mathcal{S}^{Sp}_{s, \sigma}(\mathbb{R}^n_{+}).$ In Sec. 6 we formulate the assumptions on the behaviour of~$\partial \Omega$ in a neighbourhood of the origin and prove \textit{attainability} of $\mathcal{S}^{Sp}_{s, \sigma}(\Omega).$ The proof is based on the construction of suitable trial function using the minimizer in $\mathbb{R}^n_{+}.$ Estimates of this minimizer and of its Stinga--Torrea extension are given in Sec. 7: at first we derive rough pointwise estimate of the minimizer and then we derive more accurate estimates, analogous to \cite[Theorem 2.1]{Dem_Naz}. Technical estimates used for the proof of \textit{attainability} in $\Omega$ are given in Secs.~8, 9.

\textbf{Notation}: $x \equiv (x', x_n)$ is a point in $\mathbb{R}^{n}$ or in $\Omega;$ $y \equiv (y',y_n)$ is a point in the half-space 
\begin{equation*}
\mathbb{R}^n_+ := \{ y \equiv (y',y_n) \in \mathbb{R}^{n} \  | \ y_n > 0\}.
\end{equation*}
We use the coordinates $X \equiv (x, t) \in \Omega \times \mathbb{R}_+$ dealing with the Stinga--Torrea extension from~$\Omega$ and the coordinates $Y \equiv (y,z) \in \mathbb{R}_+^n \times \mathbb{R}_+$ dealing with the extension from $\mathbb{R}^n_+.$

$\mathbb{B}_r(x)$ and $\mathbb{S}_r(x)$ are the sphere and the ball of radius $r$ centered in $x$ respectively. Also we denote $\mathbb{B}_r := \mathbb{B}_r(\mathbb{O}_{n}),$ $\mathbb{S}_r := \mathbb{S}_r(\mathbb{O}_{n}),$ $\mathbb{B}^{+}_r := \mathbb{B}_r \cap \mathbb{R}^n_{+}$ и $\mathbb{K}^{+}_r := \mathbb{B}^{+}_{2r} \setminus \mathbb{B}^{+}_{r}$ ($\mathbb{O}_{n}$ stands for the origin in $\mathbb{R}^{n}$).

Let $\varphi_{r}(y)$ be a smooth cut-off function:
\begin{equation}
\label{phi_def}
\varphi_{r}(y) := 
\begin{cases}
1,& |y|<\frac{r}{2}\\ 
0,& |y|>r
\end{cases},
\quad
|\nabla_{y} \varphi_{r}(y)| \leq \frac{C}{r}.
\end{equation}
We use letter $C$ to denote various positive constants dependent on $n,$ $s,$ $\sigma$ only. To indicate that $C$ depends on some other parameters, we write $C(\dots).$ We also write  $o_{\varepsilon}(1)$ to indicate a quantity that tends to zero as $\varepsilon \to 0.$

We use the notation $\widetilde{u}(y)$ and $\widetilde{w}(Y)$ for the odd reflections of $u(y)$ and its Stinga--Torrea extension $w(Y):$
\begin{equation*}
\widetilde{u}(y) := \begin{cases}
u(y', y_n),& y_n \geq 0\\ 
-u(y', -y_n),& y_n \leq 0
\end{cases},
\quad
\widetilde{w}(Y) := \begin{cases}
w(y', y_n, z),& y_n \geq 0\\ 
-w(y', -y_n, z),& y_n \leq 0.
\end{cases}
\end{equation*}
\section{Preliminaries}
Recall (see, for instance, \cite[Secs. 2.3.3, 4.3.2]{Triebel}) that the Sobolev spaces $H^s(\mathbb{R}^n)$ and $\widetilde{H}^s(\Omega)$ are defined  via the Fourier transform $\mathcal{F}u(\xi) := \frac{1}{(2\pi)^{n/2}}\int\limits_{\mathbb{R}^n} e^{-i\xi\cdot x}u(x)dx$:
\begin{gather*}
H^s(\mathbb{R}^n) = \left\{ u\in L_2(\mathbb{R}^n) \ \left| \   \|u\|^2_{H^s(\mathbb{R}^n)} := \int\limits_{\mathbb{R}^n} (1+|\xi|^{2s})|\mathcal{F}u(\xi)|^2d\xi  < +\infty \right. \right\};\\
\widetilde{H}^s(\Omega) = \left\{ u\in H^s(\mathbb{R}^n) \mid \supp(u) \subset \overline{\Omega} \right\}.
\end{gather*}
The fractional Laplacian $(-\Delta)^s$ in $\mathbb{R}^n$ of a function $u \in \mathcal{C}^{\infty}_0(\mathbb{R}^n)$ is defined by the identity
\begin{equation}
(-\Delta)^s u = \mathcal{F}^{-1}(|\xi|^{2s}\mathcal{F}u(\xi)), \quad 
\langle (-\Delta)^s u,u \rangle = \int\limits_{\mathbb{R}^n} |\xi|^{2s}|\mathcal{F}u(\xi)|^2d\xi.
\label{seminorm_D}
\end{equation}
Quadratic form (\ref{seminorm_D}) is well-defined on $H^s(\mathbb{R}^n),$ thus the fractional Laplacian in $\mathbb{R}^n$ can be considered as a self-adjoint operator with the quadratic form (\ref{seminorm_D}) on $H^s(\mathbb{R}^n).$

The \textit{Spectral Dirichlet Laplacian} $(-\Delta)_{Sp}^s$ is the $s$-th power of the conventional Dirichlet Laplacian in the sense of spectral theory. Its quadratic form in $\mathbb{R}^n$ coincides with (\ref{seminorm_D}), i.e. $(-\Delta)_{Sp}^s \equiv (-\Delta)^s$ in $\mathbb{R}^n.$ In case of $\Omega =\mathbb{R}^n_{+}$ the quadratic form is equal to
\begin{equation*}
\langle (-\Delta)_{Sp}^s u,u \rangle := \int\limits_{\mathbb{R}^n_{+}} |\xi|^{2s}|\widehat{\mathcal{F}}u(\xi)|^2d\xi \quad \mbox{with} \quad  \widehat{\mathcal{F}}u(\xi) := \frac{2}{(2\pi)^{n/2}}\int\limits_{\mathbb{R}^n} u(x) e^{-i\xi'\cdot x'} \sin(x_n \xi_n)dx;
\end{equation*}
for the bounded domain $\Omega$
\begin{equation}
\label{seminorm_Sp}
\langle(-\Delta)_{Sp}^s u,u\rangle := \sum\limits_{j=1}^{\infty} \lambda_j^s \langle u,\phi_j \rangle^2,
\end{equation}
here $\lambda_j$ и $\phi_j$ are eigenvalues and eigenfunctions (orthonormalized in $L_2(\Omega)$) of the Dirichlet Laplacian in $\Omega,$ respectively.

\begin{prop}[\textbf{\cite[Theorem 2]{MN2}}]
\label{DN_ineq_prop}
Let $s \in (0,1).$ Then for $u(x) \in \widetilde{H}^{s}(\Omega)$ the following inequality holds:
\begin{equation}
\label{DN_ineq}
 \langle(-\Delta)_{Sp}^s u,u\rangle \geq  \langle(-\Delta)^s u,u\rangle.
\end{equation}
If $u \not\equiv 0,$ then (\ref{DN_ineq}) holds with strict sign.
\end{prop}
Inequality (\ref{in1}) for $s \in (0,1)$ (or even for $s \in \left(0, \tfrac{n}{2}\right)$) follows from (\ref{DN_ineq}) and the general theorem by V.P. Il'in \cite[Theorem 1.2, (22)]{Ilin} about estimates of integral operators in weighted Lebesgue spaces.
 
Let $\Omega = \mathbb{R}^n,$ then for $\sigma=0$ inequality (\ref{in1}) reduces to the fractional Hardy inequality
\begin{equation}
\label{fHardy}
\langle (-\Delta)^s u, u \rangle \geq \mathcal{S}_{s,0}  \||x|^{-s}u\|^2_{L_2(\mathbb{R}^n)},
\end{equation}
and for $\sigma = s$ it reduces to the fractional Sobolev inequality
\begin{equation}
\label{fSob}
\langle (-\Delta)^s u, u \rangle \geq \mathcal{S}_{s,s} \|u\|^2_{L_{2^*_{s}}(\mathbb{R}^n)}.
\end{equation}
The explicit values of $\mathcal{S}_{s,0}$ and $\mathcal{S}_{s,s}$ have been computed in \cite{Herbst} and \cite{Tav} respectively. The explicit value of $\mathcal{S}_{s,\sigma}(\mathbb{R}^n)$ for arbitrary $\sigma \in (0,1)$ is unknown.

Using (\ref{fSob}), we define $\mathcal{D}^s(\mathbb{R}^n)$ and $\widetilde{\mathcal{D}}^s(\Omega)$ spaces with the $\langle(-\Delta)_{Sp}^s u,u\rangle$ norm:
\begin{gather*}
\mathcal{D}^s(\mathbb{R}^n) := \left\{ u\in L_{2_s^*}(\mathbb{R}^n) \ \left| \  \langle (-\Delta)^s u, u \rangle < \infty \right. \right\};\\
\widetilde{\mathcal{D}}^s(\Omega) := \left\{ u \in \mathcal{D}^s(\mathbb{R}^n) \mid u \equiv 0 \ \mbox{outside of} \ \overline{\Omega} \right\}.
\end{gather*}
The space $\widetilde{\mathcal{D}}^s(\Omega)$ also can be defined as the closure of $\mathcal{C}_0^{\infty}(\Omega)$ with respect to the norm, generated by $\langle(-\Delta)_{Sp}^s u,u\rangle.$ Obviously $\widetilde{\mathcal{D}}^s(\mathbb{R}^n_+) \cap L_2(\mathbb{R}^n_+) = \widetilde{H}^{s}(\mathbb{R}^n_+),$ and for any bounded $\Omega$ the Friedrichs inequality provides $\widetilde{\mathcal{D}}^s(\Omega) \equiv \widetilde{H}^{s}(\Omega).$ 

We recall that the \textit{Spectral Dirichlet Laplacian} $(-\Delta)_{Sp}^s$ can be derived via the Stinga--Torrea extension \cite{Stinga}: the Dirichlet problem 
\begin{equation}
\label{eq_st_tor}
\mathcal{L}_s [w](X) \equiv -div(t^{1-2s}\nabla_{X} w(x,t)) = 0 \quad \mbox{in} \quad  \Omega \times \mathbb{R}_{+}; \quad \left.w\right|_{t=0} = u; \quad  \left.w\right|_{x \in \partial\Omega} = 0
\end{equation}
has a unique solution $w_{sp}$ with finite energy 
\begin{equation}
\label{Navier_ext_functional}
\mathcal{E}_s[w] := \int\limits_0^{+\infty} \int\limits_{\Omega}t^{1-2s}|\nabla_{X} w(x,t)|^2dxdt.
\end{equation}
In addition, the following relation holds in the sense of distributions:
\begin{equation}
\label{bound_w}
(-\Delta)_{Sp}^s u(x) 
=
C_s \frac{\partial w_{sp}}{\partial \nu_s} (x,0) 
:= 
-C_s \lim_{t \to 0_{+}} t^{1-2s} \partial_{t}w_{sp}(x,t)
\quad \mbox{with} \quad C_s := \tfrac{4^s\Gamma(1+s)}{2s \cdot \Gamma(1-s)}.
\end{equation}
Moreover, $w_{sp}$ is the minimizer of (\ref{Navier_ext_functional}) over the space
\begin{equation*}
\mathfrak{W}_s(\Omega) := \left\{ w(X) \mid \mathcal{E}_s[w]<+\infty, \left.w\right|_{t=0} = u, \left.w\right|_{x \in \partial\Omega} = 0 \right\},
\end{equation*}
and the quadratic form (\ref{seminorm_Sp}) can be expressed in terms of $\mathcal{E}_s[w_{sp}]$ (see, e.g., \cite[(2.6)]{MN3}):
\begin{equation}
\label{norm_u_w}
\langle(-\Delta)_{Sp}^s u,u\rangle = C_s \mathcal{E}_s \left[w_{sp}\right].
\end{equation}
We call any function $w(X) \in \mathfrak{W}_s(\Omega)$ an \textbf{admissible} extension of $u(x).$ Obviously, for any admissible extension $w$ we have $\mathcal{E}_s \left[w\right] \geq \mathcal{E}_s \left[w_{sp}\right].$  As we noted above, for $\Omega = \mathbb{R}^n$ the \textit{Spectral Dirichlet Laplacian} coincides with the fractional Laplacian $(-\Delta)^s$ in $\mathbb{R}^n,$ and its extension (the Caffarelli--Silvestre extension) was introduced earlier in \cite{Caffarelli}.

Attainability of $\mathcal{S}_{s, \sigma}(\Omega)$ is equivalent to the existence of minimizer for the functional~$\mathcal{I}_{\sigma, \Omega}$:
\begin{equation}
\label{fu1}
\mathcal{I}_{\sigma, \Omega}[u] := \frac{\langle(-\Delta)_{Sp}^s u,u\rangle}{\| |x|^{\sigma-s} u  \|^2_{L_{2^*_{\sigma}}(\Omega)}}.
\end{equation}
Standard variational argument shows that each minimizer of (\ref{fu1}) solves the following problem (up to the multiplication by a constant)
\begin{equation}
\label{main_equation}
(-\Delta)_{Sp}^s u(x) = \frac{|u|^{2^*_{\sigma} - 2}u(x)}{|x|^{(s-\sigma)2^*_{\sigma}}} \quad \mbox{in} \quad \Omega, \quad u\in \widetilde{\mathcal{D}}^s(\Omega).
\end{equation}

The \textbf{$s$-Kelvin transform} in $\mathfrak{W}(\mathbb{R}^n)$ is defined by formula
\begin{equation}
\label{s_kelvin}
w^{*}(X) := \frac{1}{|X|^{n-2s}} w\Bigl(\frac{X}{|X|^2}\Bigr) \quad \forall X \equiv (x,t) \in\mathbb{R}^n \times \mathbb{R}_{+} \setminus \{\mathbb{O}_{n+1}\}.
\end{equation}
The following properties are hold for the $s$-Kelvin transform (see, e.g., \cite[Proposition 2.6]{Fall}):
\begin{equation*}
\left\{
\begin{array}{ll}
\mathcal{L}_s \left[w^*\right](X) = |X|^{-n-2s-2}  \mathcal{L}_s \left[w \right] \bigl(\frac{X}{|X|^2}\bigr), \quad  & \forall X  \equiv (x,t) \in \mathbb{R}^n \times \mathbb{R}_{+} \setminus \{ \mathbb{O}_{n+1}\}, \\
\frac{\partial w^*}{\partial \nu_s} (x,0) \equiv |x|^{-n-2s} \frac{\partial w}{\partial \nu_s} \bigl(\frac{x}{|x|^2},0 \bigr),  & \forall x \in \mathbb{R}^n \setminus \{\mathbb{O}_{n}\}.
\end{array}
\right. 
\end{equation*}
The relation 
\begin{equation*}
\frac{\partial w^*}{\partial \nu_s} (x,0) 
\equiv 
|x|^{-n-2s} \frac{\partial w}{\partial \nu_s} \Bigl(\frac{x}{|x|^2},0 \Bigr) 
= 
\frac{w^{2^*_{\sigma} - 1} \bigl(\frac{x}{|x|^2},0 \bigr) |x|^{-n-2s}}{ \bigl|\frac{x}{|x|^2}\bigr|^{(s-\sigma)2^*_{\sigma}}} 
=
\frac {\left(w^*\right)^{2^*_{\sigma}-1}(x,0)}{|x|^{2^*_{\sigma}(s-\sigma)}},
\end{equation*}
shows that the problem (\ref{main_equation}) is invariant under the $s$-Kelvin transform. This fact allows us to derive estimates of $w$ near the origin and at infinity from each other.

In what follows, we need the following propositions:
\begin{prop}[\textbf{\cite[Theorem 3]{MN4}}]
\label{prop_musina}
Let $u(x) \in \widetilde{\mathcal{D}}^{s}(\Omega),$ $s \in (0,1).$ Then $|u(x)| \in \widetilde{\mathcal{D}}^{s}(\Omega)$ and 
\begin{equation*}
 \langle(-\Delta)_{Sp}^s u,u\rangle \geq  \langle(-\Delta)_{Sp}^s |u|, |u|\rangle.
\end{equation*}
Moreover, if both the positive and the negative parts of $u$ are non-trivial, then strict inequality holds.
\end{prop}
\noindent The proof in \cite{MN4} is given for bounded domains, but works for unbounded domains without any changes.
\begin{prop}[\textbf{\cite[Lemma 2.6]{Capella}}, \textbf{\cite[Proposition A.1]{MN6}}]
\label{prop_max}
Let be $s \in (0,1),$ $u \not\equiv 0,$ $u(x) \in \widetilde{\mathcal{D}}^{s}(\Omega)$ or $u(x) \in \mathcal{D}^{s}(\mathbb{R}^n),$ and $(-\Delta)_{Sp}^s u \geq 0$ holds in the sense of distributions. Then $u > 0$ for any compact $K \subset \Omega$ (or $K \subset \mathbb{R}^n$ respectively).
\end{prop}
\noindent Proposition \ref{prop_musina} shows that the substitution $u \to |u|$ decreases $\mathcal{I}_{\sigma, \Omega}.$ Therefore, if $u$ is a minimizer of (\ref{fu1}), then the right-hand side of (\ref{main_equation}) is non-negative. Thus the maximum principle from the Proposition \ref{prop_max} shows that $u$ preserves a sign.
\begin{prop}[\textbf{\cite[Proposition 3]{MN4}}]
\label{prop_limit}
Let $u(x) \in \widetilde{\mathcal{D}}^s(\Omega)$ and $u_{\rho}(x) := \rho^{\frac{n-2s}{2}} u(\rho x).$ Then
\begin{equation*}
\langle (-\Delta)^s u, u \rangle = \lim_{\rho \to \infty} \langle (-\Delta)_{\Omega, Sp}^s u_{\rho}, u_{\rho} \rangle.
\end{equation*}
\end{prop}

\section{Non-existence results}
In this section we consider the case $\mathbb{O}_n \in \Omega.$
\begin{theorem}
\label{th_non_exist}
Let $\mathbb{O}_n \in \Omega$ and $\widetilde{\mathcal{D}}^s(\Omega) \neq \mathcal{D}^{s}(\mathbb{R}^n).$
\begin{enumerate}
\item The constant $\mathcal{S}^{Sp}_{s, \sigma}(\Omega)$ is not attained.
\item If $\Omega$ is star-shaped with respect to $\mathbb{O}_n,$ then the only non-negative solution of (\ref{main_equation}) is~$u \equiv 0.$
\end{enumerate}
\end{theorem}
\begin{proof}
1. In the local case $s = 1$ this statement is well-known. We adapt it for the non-local case. At first we notice that $\mathcal{S}_{s, \sigma}(\mathbb{R}^n)$ can be approximated with $\mathcal{C}_0^{\infty}(\mathbb{R}^n)$ functions since such functions are dense in $\mathcal{D}^{s}(\mathbb{R}^n).$ Proposition \ref{prop_limit} shows that for each $u \in \mathcal{C}_0^{\infty}(\mathbb{R}^n)$ we have $u_{\rho} \in \widetilde{\mathcal{D}}^{s}(\Omega)$ for sufficiently large~$\rho$ and the following relation holds:
\begin{equation*}
\lim\limits_{\rho \to \infty} \mathcal{I}_{\sigma, \Omega}[u_{\rho}] = \mathcal{I}_{\sigma, \mathbb{R}^n}[u].
\end{equation*}
This means that $\mathcal{S}^{Sp}_{s, \sigma}(\Omega) \leq \mathcal{S}_{s, \sigma}(\mathbb{R}^n).$

Let $\mathcal{S}^{Sp}_{s, \sigma}(\Omega)$ be attained on some $u \in \widetilde{\mathcal{D}}^{s}(\Omega).$ We extend $u$ by zero to derive some minimizer in~$\mathbb{R}^n:$ inequality (\ref{DN_ineq}) gives
\begin{equation*}
\mathcal{I}_{\sigma, \mathbb{R}^n}[u] \leq \mathcal{S}^{Sp}_{s, \sigma}(\Omega) \leq \mathcal{S}_{s, \sigma}(\mathbb{R}^n)
\end{equation*}
what is impossible due to $\widetilde{\mathcal{D}}^s(\Omega) \neq \mathcal{D}^{s}(\mathbb{R}^n)$ and the maximum principle from the Proposition~\ref{prop_max}.

2. To prove statement we invent a non-local variant of the Pohozhaev identity for $(-\Delta)_{Sp}^s$ (see~\cite{Ros_Oton} for $(-\Delta)^s$ in $\mathbb{R}^n$). Note that each solution of (\ref{main_equation}) has a singularity at the origin, but is smooth outside the neighborhood of the origin. Integrating by parts, we derive from~(\ref{eq_st_tor}) (here $\eta_{\varepsilon}(x) := 1- \varphi_{\varepsilon}(x),$ where $\varphi_{\varepsilon}(x)$ was introduced in (\ref{phi_def})):
\begin{multline*}
0
= 
C_s \int\limits^{+\infty}_{0} \int \limits_{\Omega} div \left( t^{1-2s} \nabla_{X} w(X) \right)  \langle X, \nabla_{X} w(X) \rangle \eta_{\varepsilon}(x) dX
\\=
\int \limits_{\Omega} \frac{u^{2^*_{\sigma} - 1}(x)}{|x|^{(s-\sigma)2^*_{\sigma}}} \langle x, \nabla_{x} u(x) \rangle \eta_{\varepsilon}(x) dx 
+
C_s \int\limits^{+\infty}_{0} \int \limits_{\partial\Omega} t^{1-2s} \langle \nabla_{x} w(X), \vec{\,n} \rangle \langle x, \nabla_{x} w(X) \rangle \eta_{\varepsilon}(x) dX
\\-
C_s \int\limits^{+\infty}_{0} \int \limits_{\Omega} t^{1-2s} |\nabla_{X} w(X)|^2 \eta_{\varepsilon}(x) dX
-
\frac{C_s}{2} \int\limits^{+\infty}_{0} \int \limits_{\Omega} t^{1-2s} \langle X, \nabla_X \left( |\nabla_{X} w(X)|^2 \right) \rangle \eta_{\varepsilon}(x) dX
\\-
C_s \int\limits^{+\infty}_{0} \int \limits_{\Omega} t^{1-2s} \langle \nabla_{x} w(X), \nabla_{x} \eta_{\varepsilon}(x) \rangle \langle X, \nabla_{X} w(X) \rangle dX =: B_1 + B_2 + B_3 + B_4 + B_5, 
\end{multline*}
$B_1$ and $B_2$ contain $\nabla_{x}$ only since $\left.w_t\right|_{x \in \partial \Omega} = 0$ and $\left.tw_t\right|_{t = 0} = 0$ due to (\ref{bound_w}). Further,
\begin{multline*}
B_1
=
\int \limits_{\Omega} \sum\limits_{i = 1}^{n} \frac{\left[u^{2^*_{\sigma}}(x)\right]_{x_i}}{2^*_{\sigma}} \frac{x_i \eta_{\varepsilon}(x)}{|x|^{(s-\sigma)2^*_{\sigma}}} \, dx
\\=
\int \limits_{\Omega} \frac{u^{2^*_{\sigma}}(x)}{2^*_{\sigma}} \sum\limits_{i = 1}^{n} \left( \frac{\eta_{\varepsilon}(x)}{|x|^{(s-\sigma)2^*_{\sigma}}} - \frac{2^*_{\sigma}(s-\sigma) \cdot x^2_i \eta_{\varepsilon}(x)}{|x|^{(s-\sigma)2^*_{\sigma} + 2}} \right) dx
+
\int \limits_{\Omega} \frac{u^{2^*_{\sigma}}(x)}{2^*_{\sigma}} \sum\limits_{i = 1}^{n}\frac{x_i \left[\eta_{\varepsilon}(x)\right]_{x_i}}{|x|^{(s-\sigma)2^*_{\sigma}}} \, dx
\\=
\left( \frac{n}{2^*_{\sigma}} -  (s-\sigma) \right) \int \limits_{\Omega} \frac{u^{2^*_{\sigma}}(x) (1- \varphi_{\varepsilon}(x))}{|x|^{(s-\sigma)2^*_{\sigma}}} \, dx 
-
\int \limits_{\Omega} \frac{u^{2^*_{\sigma}}(x)}{2^*_{\sigma}} \sum\limits_{i = 1}^{n}\frac{x_i \left[\varphi_{\varepsilon}(x)\right]_{x_i}}{|x|^{(s-\sigma)2^*_{\sigma}}} \, dx
\\=
-\frac{n-2s}{2} \langle(-\Delta)_{Sp}^s u, u \rangle
-
\frac{n-2s}{2} \int \limits_{\Omega} \frac{u^{2^*_{\sigma}}(x) \varphi_{\varepsilon}(x)}{|x|^{(s-\sigma)2^*_{\sigma}}} \, dx 
-
\int \limits_{\Omega} \frac{u^{2^*_{\sigma}}(x)}{2^*_{\sigma}} \sum\limits_{i = 1}^{n}\frac{x_i \left[\varphi_{\varepsilon}(x)\right]_{x_i}}{|x|^{(s-\sigma)2^*_{\sigma}}} \, dx.
\end{multline*}
Since $\left.w\right|_{x \in \partial\Omega} = 0,$ vectors $\nabla_{x} w(X)$ and $\vec{\,n}$ are collinear, what gives
\begin{equation*}
B_2 
=
C_s \int\limits^{+\infty}_{0} \int \limits_{\partial\Omega} t^{1-2s} \langle x, \vec{\,n} \rangle \cdot |\nabla_{x} w(X)|^2 \eta_{\varepsilon}(x) dX.
\end{equation*}
For $B_3$ we have 
\begin{equation*}
B_3 = -\langle(-\Delta)_{Sp}^s u, u \rangle + C_s \int\limits^{+\infty}_{0} \int \limits_{\Omega} t^{1-2s} |\nabla_{X} w(X)|^2 \varphi_{\varepsilon}(x) dX.
\end{equation*}
Integrating by parts in $B_4$ we obtain (using $\left.t^{2 - 2s} |w_{t}(X)|^2\right|_{t = 0} = 0$):
\begin{multline*}
B_4 
=
- \frac{C_s}{2} \int\limits^{+\infty}_{0} \int \limits_{\partial \Omega} t^{1-2s} \langle x, \vec{\,n} \rangle |\nabla_{x} w(X)|^2 \eta_{\varepsilon}(x) dX 
\\+
 \frac{C_s (n-2s+2)}{2} \int\limits^{+\infty}_{0} \int \limits_{\Omega}  t^{1-2s} |\nabla_{X} w(X)|^2  \left( \eta_{\varepsilon}(x) + \langle x, \nabla_x \eta_{\varepsilon}(x) \rangle \right) dX
\\=
- \frac{C_s}{2} \int\limits^{+\infty}_{0} \int \limits_{\partial \Omega} t^{1-2s} \langle x, \vec{\,n} \rangle |\nabla_{x} w(X)|^2 \eta_{\varepsilon}(x) \, dX 
+
\frac{n-2s+2}{2} \langle(-\Delta)_{Sp}^s u, u \rangle
\\-
\frac{C_s}{2} \int\limits^{+\infty}_{0} \int \limits_{\Omega}  t^{1-2s} |\nabla_{X} w(X)|^2  \left( (n-2s+2)\varphi_{\varepsilon}(x) + \langle x, \nabla_x \varphi_{\varepsilon}(x) \rangle \right) dX.
\end{multline*}
Summing up, we get
\begin{multline*}
\frac{C_s}{2} \int\limits^{+\infty}_{0} \int \limits_{\partial \Omega} t^{1-2s} \langle x, \vec{\,n} \rangle |\nabla_{x} w(X)|^2 dX
=
\frac{C_s}{2} \int\limits^{+\infty}_{0} \int \limits_{\partial \Omega} t^{1-2s} \langle x, \vec{\,n} \rangle |\nabla_{x} w(X)|^2 \varphi_{\varepsilon}(x) \, dX
\\+
\frac{n-2s}{2} \int \limits_{\Omega} \frac{u^{2^*_{\sigma}}(x) \varphi_{\varepsilon}(x)}{|x|^{(s-\sigma)2^*_{\sigma}}} \, dx 
+
\int \limits_{\Omega} \frac{u^{2^*_{\sigma}}(x)}{2^*_{\sigma}} \sum\limits_{i = 1}^{n}\frac{x_i \left[\varphi_{\varepsilon}(x)\right]_{x_i}}{|x|^{(s-\sigma)2^*_{\sigma}}} \, dx
\\+
\frac{C_s}{2} \int\limits^{+\infty}_{0} \int \limits_{\Omega}  t^{1-2s} |\nabla_{X} w(X)|^2  \left( (n-2s)\varphi_{\varepsilon}(x) + \langle x, \nabla_x \varphi_{\varepsilon}(x) \rangle \right) dX
\\-
C_s \int\limits^{+\infty}_{0} \int \limits_{\Omega} t^{1-2s} \langle \nabla_{x} w(X), \nabla_{x} \varphi_{\varepsilon}(x) \rangle \langle X, \nabla_{X} w(X) \rangle \, dX.
\end{multline*}
The right-hand side of this equality tends to zero as $\varepsilon \to 0,$ therefore, the left-hand side is zero. The assumption that $\Omega$ is star-shaped around $\mathbb{O}_n$ gives $\langle x, \vec{\,n} \rangle > 0,$ thus $\nabla_x w = 0$ on~$\partial \Omega.$ Integrating by parts, we get
\begin{equation*}
0 = \int\limits^{+\infty}_{0} \int \limits_{\Omega} div(t^{1-2s} \nabla_{X}w(X)) dX = 
\int \limits_{\Omega}   \frac{u^{2^*_{\sigma} -1}(x)}{|x|^{(s-\sigma)2^*_{\sigma}}} \, dX + \lim \limits_{t \to \infty} \int \limits_{\Omega}   div(t^{1-2s} \nabla_{X}w(X)) \, dx.
\end{equation*}
The second term in the right-hand side is zero (for more details see (\ref{fourier_rep_omega}) and (\ref{bes_est_omega}) in Sec. 6), but the first term is not equal to zero for $u(x) \geq 0, u \not \equiv 0,$ what is impossible.\qedhere
\end{proof}
\begin{remark}
The second statement of the Theorem \ref{th_non_exist} is also valid in the case $0 \in \partial \Omega.$
\end{remark}
Below we assume that $0 \in \partial \Omega.$ We also can assume that $\||x|^{\sigma-s} u\|_{L_{2^*_{\sigma}}(\Omega)} = 1$ with $u(x) > 0$ due to the invariance of (\ref{fu1}) under dilations and multiplications by a constant.
\section{Estimates of Green functions}
The simplest problem involving the fractional Laplacian in $\mathbb{R}^n_{+}$ is:
\begin{equation}
(-\Delta)_{Sp}^s u(y) = h(y) \quad \mbox{in} \quad \mathbb{R}^n_{+}.
\label{neumann_u_h_Rn_+}
\end{equation}
The boundary value problem (BVP) (\ref{eq_st_tor}) in $\mathbb{R}^n_{+}$ looks like 
\begin{equation}
\mathcal{L}_s [w] \left(Y\right) \equiv -div(z^{1-2s}\nabla_{Y} w\left(Y\right)) = 0 \quad \mbox{in} \quad  \mathbb{R}^n_{+} \times \mathbb{R}_{+}; \quad \left.w\right|_{z=0} = u; \quad  \left.w\right|_{y_n = 0} = 0.
\label{neumann_w_u_Rn_+}
\end{equation}
The Stinga--Torrea extension $w\left(Y\right)$ can be derived from $h(y)$ by solving of the BVP
\begin{equation}
\mathcal{L}_s [w]\left(Y\right) = 0 \quad \mbox{in} \quad  \mathbb{R}^n_{+} \times \mathbb{R}_{+}; \quad C_s \frac{\partial w}{\partial \nu_s} (y,0) = h(y); \quad  \left.w\right|_{y_n = 0} = 0.
\label{neumann_w_h_Rn_+}
\end{equation}
\begin{lemma}
The Green functons of problems (\ref{neumann_u_h_Rn_+})-(\ref{neumann_w_h_Rn_+}) are as follows:
\begin{align}
\mbox{For (\ref{neumann_w_h_Rn_+}):}& \quad G_s(Y, \xi) := \frac{\widetilde{C}_{n,s}}{\left(|y - \xi|^2 + z^2\right)^{\frac{n-2s}{2}}} \left(1 - \left[1 + \frac{4y_n\xi_n}{|y - \xi|^2 + z^2}\right]^{\frac{2s-n}{2}}\right); \label{green_w_h_Rn+} \\
\mbox{For (\ref{neumann_w_u_Rn_+}):}& \quad {\Gamma}_s(Y, \xi) := \frac{\widehat{C}_{n,s} z^{2s}}{\left(|y - \xi|^2 + z^2\right)^{\frac{n+2s}{2}}} \left(1 - \left[1 + \frac{4y_n\xi_n}{|y - \xi|^2 + z^2}\right]^{-\frac{n+2s}{2}}\right); \label{green_w_u_Rn+} \\
\mbox{For (\ref{neumann_u_h_Rn_+}):}& \quad G_s(y, \xi) := G_s(y, 0, \xi). \label{green_u_h_Rn+}
\end{align}
\end{lemma}
\begin{proof}
To obtain required Green functions, we consider the odd reflections $\widetilde{u}(y)$ and $\widetilde{w}\left(Y\right):$ $\widetilde{w}\left(Y\right)$ is the Stinga--Torrea extension of $\widetilde{u}(y)$ due to $\left.w\right|_{y_n = 0} = 0.$ In \cite{Caffarelli} and \cite[Remark~3.10]{Cabre} the Green functions in $\mathbb{R}^n$ were calculated for two problems: for the BVP
\begin{equation*}
-div(t^{1-2s}\nabla_{X} \widetilde{w}(X)) = 0 \quad \mbox{in} \quad  \mathbb{R}^n \times \mathbb{R}_{+}; \quad C_s \frac{\partial \widetilde{w}}{\partial \nu_s} (x,0) = \widetilde{h}(x)
\end{equation*}
we have the Green function $\widetilde{G}_s(X):$
\begin{equation}
\label{green_function_g}
\widetilde{w}(X) = \int\limits_{\mathbb{R}^n} \widetilde{G}_s(x - \xi, t) \widetilde{h}(\xi) d\xi \quad \mbox{with} \quad \widetilde{G}_s(X) :=  \frac{\widetilde{C}_{n,s}}{\left(x^2 + t^2\right)^{\frac{n-2s}{2}}};
\end{equation}
for the BVP 
\begin{equation*}
-div(t^{1-2s}\nabla_{X} \widetilde{w}(X)) = 0 \quad \mbox{in} \quad  \mathbb{R}^n \times \mathbb{R}_{+}; \quad \left.\widetilde{w}\right|_{t=0} = \widetilde{u}
\end{equation*}
we have the Green function $\widetilde{\Gamma}_s(X):$
\begin{equation}
\label{green_function_gamma}
\widetilde{w}(X) = \int\limits_{\mathbb{R}^n} \widetilde{\Gamma}_s(x - \xi, t) \widetilde{u}(\xi) d\xi \quad \mbox{with} \quad \widetilde{\Gamma}_s(X) :=  \frac{\widehat{C}_{n,s} t^{2s}}{\left(x^2 + t^2\right)^{\frac{n+2s}{2}}}.
\end{equation}
Required representation (\ref{green_w_h_Rn+}) follows from (\ref{green_function_g}) and from the identity
\begin{equation*}
G_s(Y, \xi) = \widetilde{G}_s(y', y_n, t, \xi) - \widetilde{G}_s(y', -y_n, t, \xi) \quad \mbox{with} \quad y_n > 0.
\end{equation*}
Similarly, (\ref{green_w_u_Rn+}) follows from (\ref{green_function_gamma}),  the representation (\ref{green_u_h_Rn+}) is obvious.
\end{proof}
\begin{lemma}
For any $\mathfrak{b} \in [0,1]$ the Green functions $G_s(Y, \xi)$ and ${\Gamma}_s(Y, \xi)$ admit the following estimates: 
\begin{equation}
\label{green_est}
G_s(Y, \xi) \leq \frac{C y^{\mathfrak{b}}_n\xi^{\mathfrak{b}}_n }{\left(|y - \xi|^2 + z^2\right)^{\frac{n-2s+2\mathfrak{b}}{2}}} \quad \mbox{and} \quad
{\Gamma}_s(Y, \xi) \leq \frac{C y^{\mathfrak{b}}_n\xi^{\mathfrak{b}}_n z^{2s}}{\left(|y - \xi|^2 + z^2\right)^{\frac{n+2s+2\mathfrak{b}}{2}}}.
\end{equation}
Also, $\nabla_{Y} G_s(Y, \xi)$ can be estimated as follows:
\begin{equation}
\label{green_grad_est}
|\nabla_{Y} G_s(Y, \xi)| \leq \frac{C}{\left(|y - \xi|^2 + z^2\right)^{\frac{n-2s+1}{2}}} \cdot \min \left(1, \frac{y_n\xi_n }{|y - \xi|^2 + z^2} + \frac{\xi_n}{\sqrt{|y - \xi|^2 + z^2}}\right).
\end{equation}
\end{lemma}
\begin{proof}
The estimate for $G_s$ follows from the interpolation of two inequalities:
\begin{equation*}
G_s(Y, \xi) \leq \frac{C}{\left(|y - \xi|^2 + z^2\right)^{\frac{n-2s}{2}}} \quad \mbox{and} \quad G_s(Y, \xi) \leq \frac{C y_n\xi_n }{\left(|y - \xi|^2 + z^2\right)^{\frac{n-2s+2}{2}}},
\end{equation*}
the first one is obvious, and the second follows from the mean value theorem:
\begin{equation}
\label{lag_est}
1 - \left[1 + \frac{4y_n\xi_n}{|y - \xi|^2 + z^2}\right]^{\frac{2s-n}{2}} \leq \frac{C y_n\xi_n}{|y - \xi|^2 + z^2}.
\end{equation}
The estimate for ${\Gamma}_s$ can be obtained in the same way.

The gradient $\nabla_{Y} G_s(Y, \xi)$ is given by the formulae (here $i \in [1:n-1]$)
\begin{equation*}
\begin{pmatrix}
 \partial_{z} G_s(Y, \xi)\\
 \partial_{y_i} G_s(Y, \xi)\\
 \partial_{y_n} G_s(Y, \xi)
\end{pmatrix}
= C \cdot
\begin{pmatrix}
\frac{z}{\left(|y - \xi|^2 + z^2\right)^{\frac{n-2s+2}{2}}} \left(1 - \left[1 + \frac{4y_n\xi_n}{|y - \xi|^2 + z^2}\right]^{\frac{2s - n-2}{2}}\right) \\
\frac{y_i - \xi_i}{\left(|y - \xi|^2 + z^2\right)^{\frac{n-2s+2}{2}}} \left(1 - \left[1 + \frac{4y_n\xi_n}{|y - \xi|^2 + z^2}\right]^{\frac{2s - n-2}{2}}\right)\\
\frac{y_n - \xi_n}{\left(|y - \xi|^2 + z^2\right)^{\frac{n-2s+2}{2}}} - \frac{y_n + \xi_n}{\left(|y' - \xi'|^2 + |y_n + \xi_n|^2 + z^2\right)^{\frac{n-2s+2}{2}}}
\end{pmatrix},
\end{equation*}
therefore, the first part of (\ref{green_grad_est}) is obvious. The second part for $\partial_{z} G_s$ и $\partial_{y_i} G_s$ can be derived using the analogue of (\ref{lag_est}). Inequality for $\partial_{y_n}G_s$ follows from the inequality (recall that $\xi_n > 0$ and $ y_n>0$) 
\begin{equation*}
|\partial_{y_n} G_s(Y, \xi)| \leq \frac{|y_n-\xi_n|}{\left(|y - \xi|^2 + z^2\right)^{\frac{n-2s+2}{2}}} \left(1 - \left[1 + \frac{4y_n\xi_n}{|y - \xi|^2 + z^2}\right]^{\frac{2s-n-2}{2}}\right) + \frac{2\xi_n}{\left(|y - \xi|^2 + z^2\right)^{\frac{n-2s+2}{2}}}
\end{equation*}
and the analogue of (\ref{lag_est}) for the expression in large brackets.
\end{proof}
\section{Attainability of $\mathcal{S}^{Sp}_{s, \sigma}(\mathbb{R}^n_{+})$}
In this section we prove the existence of the minimizer for the functional (\ref{fu1}) in the case  $\Omega = \mathbb{R}^n_+$ and discuss its properties.
\begin{theorem}
\label{Rn+_existence}
For $\Omega = \mathbb{R}^n_+$ there exists a minimizer of the functional (\ref{fu1}).
\end{theorem}
\begin{proof}
We follow the scheme given in \cite[Theorem 3.1]{Naz_Cone} and based on the concentration-compactness principle of Lions \cite{Lions}. Consider a minimizing sequence $\{ u_k \}$ for (\ref{fu1}). As was mentioned in Sec. 2, we can assume that $u_k(y) \geq 0$ и $\| |y|^{\sigma-s} u_k \|_{L_{2^*_{\sigma}}\left(\mathbb{R}^n_+\right)} = 1.$ We also denote the Stinga--Torrea extensions as $w_k\left(Y\right)$ and define functions $U_k(y)$ as
\begin{equation}
\label{U_k_def}
U_k(y) :=  \int\limits_{0}^{+\infty} z^{1-2s} |\nabla_{Y} w_k(Y)|^2 dz.
\end{equation}
Since $\{ u_k \}$ is bounded in $\widetilde{\mathcal{D}}^s(\mathbb{R}^n_+),$ $w_k\left(Y\right)$ are uniformly bounded in $\mathfrak{W}_s(\mathbb{R}^n_{+})$ as well as $U_k$ и $||y|^{\sigma-s} u_k|^{2^*_{\sigma}}$ are uniformly bounded in $L_1(\mathbb{R}^n_+).$ Without loss of generality, we assume that:
\begin{itemize}
\item $u_k \rightharpoondown u$ in $\widetilde{\mathcal{D}}^s(\mathbb{R}^n_+);$
\item $\nabla_{Y} w_k \rightharpoondown \nabla_{Y} w$ in $L_2(\mathbb{R}^n_+ \times \mathbb{R}_+, z^{1-2s}),$ and $w$ is an admissible extension of $u;$
\item $||y|^{\sigma-s} u_k|^{2^*_{\sigma}}$ weakly converges to a measure $\mu$ on $\overline{\mathbb{R}^n_+};$
\item $U_k$ weakly converges to a measure  $\mathcal{M}$ on $\overline{\mathbb{R}^n_+};$
\end{itemize}
where $\overline{\mathbb{R}^n_+}$ is a one-point compactification of $\mathbb{R}^n_+.$

Embedding $\widetilde{\mathcal{D}}^s_{loc}(\mathbb{R}^n_+) \hookrightarrow L_{2^*_{\sigma}, loc}(\mathbb{R}^n_+ \backslash \{\mathbb{O}_n\})$ is compact due to $2^*_{\sigma} < 2^*_{s},$ therefore $|y|^{\sigma - s} u_k \to |y|^{\sigma - s} u$ in $L_{2^*_{\sigma}, loc}(\mathbb{R}^n_+ \backslash \{\mathbb{O}_n\})$ and we have the following representation for $\mu:$
\begin{equation*}
\mu = ||y|^{\sigma - s} u|^{2^*_{\sigma}} + \alpha_0 \bm{\delta_0}(y) + \alpha_{\infty} \bm{\delta_{\infty}}(y), \quad\alpha_{0}, \alpha_{\infty} \geq 0,
\end{equation*}
here $\bm{\delta_0}(y)$ and  $\bm{\delta_{\infty}}(y)$ are Dirac delta functions at the origin and at infinity respectively. 

Our next goal is to show that the measure $\mathcal{M}$ admits the estimate:
\begin{equation}
\label{M_ineq}
\mathcal{M} \geq U + \mathcal{S}^{Sp}_{s, \sigma}(\mathbb{R}^n_+)  \alpha_0^{\tfrac{2}{2^*_{\sigma}}}  \bm{\delta_0}(y) + \mathcal{S}^{Sp}_{s, \sigma}(\mathbb{R}^n_+)  \alpha_{\infty}^{\tfrac{2}{2^*_{\sigma}}}  \bm{\delta_{\infty}}(y).
\end{equation}
Obviously, it suffices to prove that $\mathcal{M}$ majorizes separately each term in the right-hand side of~(\ref{M_ineq}) . The first estimate $\mathcal{M} \geq U$ follows from the weak convergence $\nabla_{Y} w_k \eta \rightharpoondown \nabla_{Y} w \eta$ in $L_2\left(\mathbb{R}^n_+ \times \mathbb{R}_+, z^{1-2s} \right)$ for any $\eta \in \mathcal{C}_0^{\infty}(\mathbb{R}^n_+)$ and from the weak lower semi-continuity of the weighted $L_2$-norm:
\begin{multline}
\label{meas_gr_U}
\int\limits_{\mathbb{R}^n_+} \eta^2(y) d \mathcal{M} \equiv \lim \limits_{k \to \infty} \int\limits_{\mathbb{R}^n_+} \eta^2(y) U_k(y) dy
= 
\lim \limits_{k \to \infty} \int\limits_{0}^{+\infty} \int\limits_{\mathbb{R}^n_+} z^{1-2s}|\nabla_{Y} w_k(Y) \cdot \eta(y)|^2 dY 
\\ \geq
\int\limits_{0}^{+\infty} \int\limits_{\mathbb{R}^n_+} z^{1-2s}|\nabla_{Y} w(Y) \cdot \eta(y)|^2 dY 
=
\int\limits_{\mathbb{R}^n_+} \eta^2(y) U(y) dy.
\end{multline}
To obtain the second estimate we use the trial function $\eta_{\varepsilon}(y) := \varphi_{2\varepsilon}(y):$ 
\begin{multline}
\label{delta_case}
\int\limits_{\mathbb{R}^n_+} \eta_{\varepsilon}^2(y) d\mathcal{M} \equiv 
\lim \limits_{k \to \infty} \int\limits_{\mathbb{R}^n_+} U_k \eta_{\varepsilon}^2(y) dy 
=
\lim \limits_{k \to \infty} \int\limits_{0}^{+\infty} \int\limits_{\mathbb{R}^n_+} z^{1-2s}|\nabla_{Y}\left[w_k(Y)\eta_{\varepsilon}(y)\right] - w_k(Y)\nabla_{y}\eta_{\varepsilon}(y)|^2 dY 
\\=
\lim \limits_{k \to \infty} \int\limits_{0}^{+\infty} \int\limits_{\mathbb{R}^n_+} \Bigl[ z^{1-2s}|\nabla_{Y}\left[w_k(Y)\eta_{\varepsilon}(y)\right]|^2
-
2 z^{1-2s} \nabla_{y} w_k(Y)\nabla_{y}\eta_{\varepsilon}(y) w_k(Y)\eta_{\varepsilon}(y)
\\+ 
z^{1-2s}|w_k(Y) \nabla_{y}\eta_{\varepsilon}(y)|^2 \Bigr] dY 
=: D_1 - D_2 + D_3.
\end{multline}
To estimate $D_1$ we use the Hardy--Sobolev inequality (\ref{in1}):
\begin{equation}
\label{D_1_est}
D_1 
\geq 
\mathcal{S}^{Sp}_{s, \sigma}(\mathbb{R}^n_+) \cdot  \lim \limits_{k \to \infty} \||y|^{\sigma - s} u_k\varphi_{2\varepsilon}\|^2_{L_{2^*_{\sigma}}(\mathbb{R}^n_+)}
\geq
\mathcal{S}^{Sp}_{s, \sigma}(\mathbb{R}^n_+)  \alpha_0^{\tfrac{2}{2^*_{\sigma}}}.
\end{equation}
To estimate $D_3$ we have to pass to the limit under the integral side:
\begin{lemma}
\label{D3_lemma}
The following equality holds:
\begin{equation}
\label{D3_eq}
D_3 = \int\limits_{0}^{+\infty} \int\limits_{\mathbb{R}^n_+} z^{1-2s}w^2(Y) |\nabla_{y}\eta_{\varepsilon}(y)|^2 dY.
\end{equation}
\end{lemma}
\begin{proof}
Let be $\delta \in (0,1),$ we split the integral into three parts:
\begin{equation*}
a_k + b_k + c_k := \left( \int\limits_{0}^{\delta} + \int\limits_{\delta}^{\frac{1}{\delta}} + \int\limits_{\frac{1}{\delta}}^{+\infty}\right) \int\limits_{\mathbb{R}^n_+} z^{1-2s} w^2_k(Y) |\nabla_{y}\eta_{\varepsilon}(y)|^2 dY.
\end{equation*}
For a fixed $\delta$ we can pass to the limit in $b_k:$ $w_k$ are uniformly bounded in $W^1_2\left(\mathbb{K}^{+}_{\varepsilon} \times \left[ \delta, \frac{1}{\delta} \right]\right),$ thus $w_k \to w$ in $L_2\left( \mathbb{K}^{+}_{\varepsilon} \times \left[ \delta, \frac{1}{\delta} \right]\right).$ To complete the proof, it suffices to show that 
\begin{equation}
\label{a_k_c_k}
a_k + c_k < C(\varepsilon) \cdot \delta^{1-s}
\end{equation}
To prove (\ref{a_k_c_k}) for $a_k$ we use the Green function (\ref{green_w_u_Rn+}):
\begin{gather*}
w_k\left(Y\right) = \int\limits_{\mathbb{R}^n_{+}} u_k(\xi) {\Gamma}_s(Y, \xi)  d\xi = \left( \  \int\limits_{|y - \xi| > 1} + \int\limits_{ |y - \xi| \leq 1} \right) u_k(\xi) {\Gamma}_s(Y, \xi)  d\xi =: w_{1k}\left(Y\right) + w_{2k}\left(Y\right), \\
a_k \leq
2\int\limits_{0}^{\delta} \int\limits_{\mathbb{R}^n_+} z^{1-2s} \left[ w^2_{1k}(Y) + w^2_{2k}(Y)\right] \cdot \left|\nabla_{y} \eta_{\varepsilon}(y)\right|^2 dY =: a_{1k} + a_{2k}.
\end{gather*}
Using (\ref{green_est}) for $\mathfrak{b} = 0,$ $|\nabla_{y} \varphi_{2\varepsilon}| \leq \frac{c}{\varepsilon}$ and the Cauchy--Bunyakovsky--Schwarz inequality we get:
\begin{multline*}
a_{1k} 
\leq 
\frac{C}{\varepsilon^2}
\int\limits_{0}^{\delta} z^{1-2s} \int\limits_{|y| < 2\varepsilon} \left( \  \int\limits_{|y-\xi| > 1} \frac{u_k(\xi) z^{2s}}{(|y - \xi|^2 + z^2)^{\frac{n+2s}{2}}} d\xi \right)^2 dY
\\ \leq
C \frac{\delta^{2+2s}}{\varepsilon^2} \int\limits_{|y| < 2\varepsilon} \left( \  \int\limits_{|y-\xi| > 1} \frac{u_k(\xi) |\xi|^{-s} |\xi|^{s}}{|y - \xi|^{n+2s}} d\xi \right)^2 dy
\leq 
C(\varepsilon) \delta^{2+2s} \||y|^{-s}u_k\|^2_{L_2(\mathbb{R}^n_+)} \int\limits_{1}^{+\infty} r^{-n-2s-1} dr.
\end{multline*}
Similarly, we get the estimate
\begin{multline*}
a_{2k} 
\leq 
\frac{C}{\varepsilon^2} \int\limits_{0}^{\delta} z^{-s} \int\limits_{|y| < 2\varepsilon} \left( \ \int\limits_{|y - \xi| \leq 1} \frac{u_k(\xi) z^{2s + \frac{1-s}{2}} d\xi}{(|y - \xi|^2 + z^2)^{\frac{n+2s}{2}}}\right)^2 dY
\leq  \\ \leq
C\frac{\delta^{1-s}}{\varepsilon^2} \int\limits_{|y| < 2\varepsilon} \left( \ \int\limits_{|y - \xi| \leq 1} \frac{u_k(\xi) d\xi}{|y - \xi|^{n - \frac{1-s}{2}}} \right)^2 dy.
\end{multline*}
We estimate the integrand from the right-hand side as follows:
\begin{multline*}
\left( \ \int\limits_{|y - \xi| \leq 1} \frac{u_k(\xi)}{|y - \xi|^{n - \frac{1-s}{2}}} d\xi \right)^2
\leq 
C \left( \ \int\limits_{|y - \xi| \leq 1} \frac{u_k(\xi) - u_k(y)}{|y - \xi|^{n - \frac{1-s}{2}}} d\xi + u_k(y) \right)^2
\\ \leq
C \left( \ \int\limits_{|y - \xi| \leq 1} \frac{u_k(\xi) - u_k(y)}{|y - \xi|^{n - \frac{1-s}{2}}} d\xi \right)^2 + Cu^2_k(y) 
\leq
C \int\limits_{|y - \xi| \leq 1} \frac{|u_k(\xi) - u_k(y)|^2}{|y - \xi|^{n + 2s}} d\xi \cdot \int\limits_{0}^1 r^{s} dr + Cu^2_k(y).
\end{multline*}
Finally, using (\ref{DN_ineq}) we get
\begin{equation*}
a_{2k} \leq C(\varepsilon) \delta^{1-s} \left( \langle(-\Delta)_{Sp}^s u_k,u_k\rangle \cdot \int\limits_{0}^{1} r^{s} dr + \||y|^{-s}u_k\|^2_{L_2(\mathbb{R}^n_+)} \right).
\end{equation*}
To prove (\ref{a_k_c_k}) for $c_k$ we use (\ref{green_est}) with $\mathfrak{b} = 1:$
\begin{multline*}
c_k 
\leq 
\frac{C}{\varepsilon^2} \int\limits_{\frac{1}{\delta}}^{+\infty} z^{1-2s} \int\limits_{|y| < 2\varepsilon} \left( \  \int\limits_{\mathbb{R}^n_{+}} \frac{u_k(\xi) z^{2s} |y_n| \xi_n}{(|y-\xi|^2 + z^2)^{\frac{n+2s + 2}{2}}} \, d\xi\right)^2 dY 
\\ \leq
C(\varepsilon) \||y|^{-s}u_k\|^2_{L_2(\mathbb{R}^n_+)} \int\limits_{\frac{1}{\delta}}^{+\infty} z^{1+2s} \left( \  \int\limits_{0}^{+\infty} \frac{r^{n+1+2s}}{(r^2 + z^2)^{n+2s+2}} \, dr\right) dz
\\ \leq
C(\varepsilon) \||y|^{-s}u_k\|^2_{L_2(\mathbb{R}^n_+)} \int\limits_{\frac{1}{\delta}}^{+\infty} z^{-1 - n} dz
=
C(\varepsilon) \delta^{n} \||y|^{-s}u_k\|^2_{L_2(\mathbb{R}^n_+)}.
\end{multline*}
Thus, the estimate (\ref{a_k_c_k}) is proved completely and we get (\ref{D3_eq}).
\end{proof}
Lemma \ref{D3_lemma} implies
\begin{equation}
\label{D3_in}
D_3 \leq \frac{C}{\varepsilon^2} \int\limits_{0}^{+\infty} \int\limits_{|y| < 2\varepsilon} z^{1-2s}|w(Y)|^2 dY.
\end{equation}
For $y_n \leq 2\varepsilon,$ using the inequality 
\begin{equation*}
|w(Y)|^2 = 
\left(\int\limits_{0}^{y_n} \frac{\partial w(y', t, z)}{\partial y_n} \, dt \right)^2
\leq  2\varepsilon
\int\limits_{0}^{2\varepsilon} \left(\frac{\partial w}{\partial y_n}\right)^2 dt,
\end{equation*}
we obtain
\begin{equation*}
D_3 
\leq 
\int\limits_{0}^{+\infty} \int\limits_{|y| < 2 \sqrt{2}\varepsilon} z^{1-2s}\left(\frac{\partial w}{\partial y_n}\right)^2 dY
= o_{\varepsilon}(1) \cdot \mathcal{E}_s[w].
\end{equation*}
To estimate $D_2$ we use the Cauchy--Bunyakovsky--Schwarz inequality
\begin{equation}
\label{D2_estimate}
|D_2|
\leq 
C \sqrt{D_1 \cdot D_3} = o_{\varepsilon}(1).
\end{equation}
To sum up, we have transformed (\ref{delta_case}) into
\begin{equation*}
\int\limits_{\mathbb{R}^n_+} \varphi_{2\varepsilon}^2 \, d\mathcal{M} \equiv 
\lim \limits_{k \to +\infty} \int\limits_{\mathbb{R}^n_+} U_k \varphi_{2\varepsilon}^2 \, dy 
\geq \mathcal{S}^{Sp}_{s, \sigma}(\mathbb{R}^n_+) \alpha_0^{\tfrac{2}{2^*_{\sigma}}} + o_{\varepsilon}(1),
\end{equation*}
what gives $\mathcal{M} \geq \mathcal{S}^{Sp}_{s, \sigma}(\mathbb{R}^n_+)  \alpha_0^{\tfrac{2}{2^*_{\sigma}}}  \bm{\delta_0}(y).$
\medskip

To derive the estimate at infinity, we put $\eta_{\varepsilon}(y) := 1 - \varphi_{\frac{2}{\varepsilon}}(y)$ and write (\ref{delta_case}) for it. The estimate for $D_1$ is similar to (\ref{D_1_est}). To estimate $D_3$ we use the analogue of Lemma \ref{D3_lemma}:
\begin{equation*}
D_3 \leq \varepsilon^2 \int\limits_{0}^{+\infty} \int\limits_{\mathbb{K}_{\frac{1}{\varepsilon}}} z^{1-2s}|w(Y)|^2 \, dY.
\end{equation*}
In spherical coordinates $(r, \theta_1, \dots ,\theta_{n-1})$ we have $w = 0$ for $\theta_{n-1} = 0,$ thus
\begin{gather*}
|w(Y)|^2 = |w(r, \theta_1, \dots ,\theta_{n-1})|^2 = 
\left(\int\limits_{0}^{\theta_{n-1}} \frac{\partial w(y', t, z)}{\partial \theta_{n-1}} \, dt \right)^2
\leq \pi
\int\limits_{0}^{\pi} \left(\frac{\partial w}{\partial \theta_{n-1}}\right)^2 dt \\
D_3 
\leq 
\pi^2 \varepsilon^2
\int\limits_{0}^{+\infty} \int\limits_{\mathbb{K}_{\frac{1}{\varepsilon}}} z^{1-2s}\left(\frac{\partial w}{\partial  \theta_{n-1}}\right)^2 dY 
\leq
\frac{4\pi^2 \varepsilon^2}{\varepsilon^2} \int\limits_{0}^{+\infty} \int\limits_{\mathbb{K}_{\frac{1}{\varepsilon}}} z^{1-2s} |\nabla_y w|^2 dY
= o_{\varepsilon}(1) \cdot \mathcal{E}_s[w].
\end{gather*}
Further arguments are similar to the estimate at the origin. Inequality (\ref{M_ineq}) is proved.

The end of the proof is rather standard. By dilations and multiplications on a suitable constant one can achieve
\begin{equation}
\label{u_k_prop}
\| |y|^{\sigma-s} u_k \|_{L_{2^*_{\sigma}}(\mathbb{B}^{+}_1)} = \| |y|^{\sigma-s} u_k \|_{L_{2^*_{\sigma}}(\mathbb{R}^n_+ \setminus \mathbb{B}^{+}_1)} = \frac{1}{2}.
\end{equation}
From (\ref{M_ineq}) and the fact that $w$ is an admissible extension of $u$ we get
\begin{multline}
\label{measure_ineq}
 \mathcal{S}^{Sp}_{s, \sigma}(\mathbb{R}^n_+) \left (\| |y|^{\sigma-s} u \|^2_{L_{2^*_{\sigma}}(\mathbb{R}^n_+)} + \alpha_0^{\tfrac{2}{2^*_{\sigma}}} + \alpha_{\infty}^{\tfrac{2}{2^*_{\sigma}}}\right)
\leq  
\langle(-\Delta)_{Sp}^s u,u\rangle + \mathcal{S}^{Sp}_{s, \sigma}(\mathbb{R}^n_+)  \alpha_0^{\tfrac{2}{2^*_{\sigma}}}  + \mathcal{S}^{Sp}_{s, \sigma}(\mathbb{R}^n_+) \alpha_{\infty}^{\tfrac{2}{2^*_{\sigma}}} 
\\ \leq 
\int\limits_{\mathbb{R}^n_+} 1 d\mathcal{M}
=
\mathcal{S}^{Sp}_{s, \sigma}(\mathbb{R}^n_+) \left( \ \int\limits_{\mathbb{R}^n_+} 1 d\mu\right)^{\tfrac{2}{2^*_{\sigma}}} = 
\mathcal{S}^{Sp}_{s, \sigma}(\mathbb{R}^n_+) \left(\| |y|^{\sigma-s} u \|^{2^*_{\sigma}}_{L_{2^*_{\sigma}}(\mathbb{R}^n_+)} + \alpha_0 + \alpha_{\infty} \right)^{\tfrac{2}{2^*_{\sigma}}},
\end{multline}
what can be true only if two of three terms from the right-hand side vanish. The relation~(\ref{u_k_prop}) keeps only the possibility $\alpha_0 = \alpha_{\infty} = 0,$ i.e. $u$ is a minimizer of (\ref{fu1}).
\end{proof}
\begin{remark}
The minimizer existence for any cone in $\mathbb{R}^n$ can be proved in a similar way. 
\end{remark}
We denote the obtained minimizer in $\mathbb{R}^n_+$ by $\Phi(y),$ and its Stinga--Torrea extension by~$\mathcal{W}(Y).$ Without loss of generality, we can assume that  $\| |y|^{\sigma-s} \Phi \|_{L_{2^*_{\sigma}}\left(\mathbb{R}^n_+\right)} = 1,$ therefore we have $\mathcal{E}_{s}\left[\mathcal{W}\right] = \mathcal{S}^{Sp}_{s, \sigma}(\mathbb{R}^n_{+}).$
\begin{lemma}
$\Phi(y)$ and $\mathcal{W}(Y)$ are radial in $y'$ and positive for $y_n > 0.$
\end{lemma}
\begin{proof}
The positivity of $\Phi(y)$ and $\mathcal{W}(Y)$ was proved at the end of Sec. 2. To prove the first part we show that a non-trivial partial Schwarz symmetrization on $y'$ (we denote the symmetrization of $u$ as $u^{\star})$ decreases (\ref{fu1}):
\begin{equation*}
\mathcal{I}_{\sigma, \Omega}[u] 
= 
\frac{\mathcal{E}_s \left[w_{sp}\right]}{\| |x|^{\sigma-s} u  \|^2_{L_{2^*_{\sigma}}(\Omega)}} 
\stackrel{*}{\geq}
\frac{\mathcal{E}_s \left[w^{\star}_{sp}\right]}{\| |x|^{\sigma-s} u  \|^2_{L_{2^*_{\sigma}}(\Omega)}}
\stackrel{**}{>} 
\frac{\mathcal{E}_s \left[w^{\star}_{sp}\right]}{\| |x|^{\sigma-s} u^{\star}  \|^2_{L_{2^*_{\sigma}}(\Omega)}}
\stackrel{***}{\geq}
\mathcal{I}_{\sigma, \Omega}[u^{\star}].
\end{equation*}
The inequality (*) is provided by the fact that $\mathcal{E}_s \left[w_{sp}\right]$ does not increase under symmetrization (see \cite[Theorem 2.31, p. 83]{Kawohl} for the Steiner symmetrization; partial Schwarz symmetrization can be achieved as the limit of Steiner symmetrizations). The inequality (**) follows from \cite[Theorem 3.4]{Lieb_Lo}. The fact that $w^{\star}_{sp}$ is an admissible extension for $u^{\star}$ gives (***).
\end{proof}
\begin{remark}
\label{un_min}
Minimizer of (\ref{fu1}) with $\| |y|^{\sigma-s} \Phi \|_{L_{2^*_{\sigma}}\left(\mathbb{R}^n_+\right)} = 1$ is not unique. Indeed, the functional (\ref{fu1}) is invariant with respect to dilations and multiplications by constant. Compositions of these transformations that keep $\| |y|^{\sigma-s} \Phi \|_{L_{2^*_{\sigma}}\left(\mathbb{R}^n_+\right)}$ norm give us multiple minimizers.
\end{remark}
For further discussion, we fix some minimizer and study its behaviour at the origin and at infinity:
\begin{lemma}
\label{min_estimates}
Minimizer $\Phi(y)$ and its Stinga--Torrea extension $\mathcal{W}(Y)$ admit the following estimates:
\begin{gather}
\label{W_common_ineq}
\Phi(y) \leq \frac{C y_n}{1 + |y|^{n-2s+2}},  \ y \in \mathbb{R}^n_{+}; \quad \mathcal{W}(Y) \leq \frac{C y_n}{1 + |Y|^{n-2s+2}},  \ Y \in \mathbb{R}^n_{+} \times \mathbb{R}_{+};\\
\label{GradW_common_ineq}
\mathcal{V}(y) := \int\limits_{0}^{+\infty} z^{1-2s}  |\nabla_{Y} \mathcal{W}(Y)|^2 dz \leq \frac{C}{1 + |y|^{2n-2s+2}}, \ y \in \mathbb{R}^n_{+},
\end{gather}
where constants $C$ depend on $n, s, \sigma$ and on the choice of the minimizer $\Phi.$
\end{lemma}
The proof of Lemma \ref{min_estimates} is given in Sec. 7.
\section{Attainability of $\mathcal{S}^{Sp}_{s, \sigma}(\Omega)$}
We assume that in a small ball $\mathbb{B}_{r_0}$ (centered at the origin) the surface $\partial \Omega$ is parametrized by the equation $x_n = F(x'),$ where $F \in \mathcal{C}^1,$ $F(\mathbb{O}_{n-1}) = 0$ and $\nabla_{x'} F(\mathbb{O}_{n-1}) = \mathbb{O}_{n-1}.$ Outside this ball $\partial \Omega$ can be arbitrary.

Following \cite{Dem_Naz} we assume that $\partial \Omega$ is \textbf{average concave at the origin}: for small $\tau > 0$
\begin{equation}
\label{f_convexity}
f(\tau) := \frac{1}{|\mathbb{S}_{\tau}^{n-2}|} \int\limits_{\mathbb{S}_{\tau}^{n-2}} F(y') \, dy' < 0.
\end{equation}
Obviously, $f \in \mathcal{C}^1$ for small $\tau.$ We also assume that $f$ is \textbf{regularly varying} at the origin with the exponent $\alpha \in [1, n-2s+3):$ for any $d > 0$
\begin{equation}
\label{f_regularity}
\lim\limits_{\tau \to 0} \frac{f(d\tau)}{f(\tau)} = d^{\alpha}.
\end{equation}
It is well known (see, e.g., \cite[Secs. 1.1, 1.2]{Seneta}) that (\ref{f_regularity}) entails $f(\tau) := - \tau^{\alpha} \psi(\tau)$ with the slowly varying function $\psi(\tau)$ (SVF). Note that for $\alpha = 1$ condition $F \in \mathcal{C}^1$ implies $\lim\limits_{\tau \to 0} \psi(\tau) = 0.$

We also introduce the functions
\begin{gather*}
f_1(\tau) := \frac{1}{|\mathbb{S}_{\tau}^{n-2}|} \int \limits_{\mathbb{S}_{\tau}^{n-2}} F^2(y') \, dy';  \quad
f_2(\tau) := \frac{1}{|\mathbb{S}_{\tau}^{n-2}|} \int \limits_{\mathbb{S}_{\tau}^{n-2}} |\nabla_{y'} F(y')|^2 \, dy'; \quad \\
f_3(\tau) :=\frac{1}{|\mathbb{S}_{\tau}^{n-2}|} \int \limits_{\mathbb{S}_{\tau}^{n-2}} |\nabla_{y'} F(y')| \, dy',
\end{gather*}
and assume that the following \textbf{condition} is fulfilled
\begin{equation}
\label{f2_f_relations}
\lim\limits_{\tau \to 0} \frac{f_2(\tau)}{f(\tau)} \tau = 0.
\end{equation}
\begin{remark}
In case of $\partial \Omega \in \mathcal{C}^2$ with negative mean curvature our assumptions (\ref{f_convexity})-(\ref{f2_f_relations}) are fulfilled with $\alpha = 2$ (see \cite[Remark 1]{Dem_Naz}). We also emphasize that these assumptions admit the absence of mean curvature ($\alpha < 2$) or its vanishing ($\alpha > 2$).
\end{remark}
\begin{remark}
It was shown in \cite[Sec. 4, (17)]{Dem_Naz} that (\ref{f2_f_relations}) implies
\begin{equation}
\label{f1_f_relations}
f_1(\tau) \leq C \tau |f(\tau)| \cdot o_{\tau}(1).
\end{equation}
\end{remark}
\begin{theorem}
\label{main_theorem}
Let $\partial \Omega$ satisfy (\ref{f_convexity})-(\ref{f2_f_relations}). Then the minimizer of (\ref{fu1}) exists,  i.e. the problem~(\ref{main_equation}) has a positive solution in $\Omega.$
\end{theorem}
\begin{proof}
The scheme of the proof is the same as in Theorem \ref{Rn+_existence}. Consider a minimizing sequence $\{ u_k \}$ for (\ref{fu1}). We denote the Stinga--Torrea extensions as $w_k\left(Y\right)$ and define functions $U_k(y)$ via (\ref{U_k_def}). As before, $U_k \in L_1(\Omega)$ and $||x|^{\sigma-s} u_k|^{2^*_{\sigma}} \in L_1(\Omega),$ and we can also assume that:
\begin{itemize}
\item $u_k \geq 0,$ $u_k \rightharpoondown u$ in $\widetilde{\mathcal{D}}^s(\Omega);$
\item $\nabla_{X} w_k \rightharpoondown \nabla_{X} w$ in $L_2(\Omega \times \mathbb{R}_{+}, t^{1-2s})$ and $w$ is an admissible extension of $u;$
\item $||x|^{\sigma-s} u_k|^{2^*_{\sigma}}$ weakly converges to a measure $\mu$ on $\overline{\Omega};$
\item $U_k$ weakly converges to a measure $\mathcal{M}$ on $\overline{\Omega}.$
\end{itemize}
In contrast to the case of $\mathbb{R}^n_+,$ for the bounded $\Omega$
\begin{equation*}
\mu = ||x|^{\sigma - s} u|^{2^*_{\sigma}} + \alpha_0 \bm{\delta_0}(x),
\end{equation*}
and we should show that 
\begin{equation}
\label{M_ineq_omega}
\mathcal{M} \geq U + \mathcal{S}^{Sp}_{s, \sigma}(\Omega)  \alpha_0^{\tfrac{2}{2^*_{\sigma}}}  \bm{\delta_0}(x).
\end{equation}
The estimate $\mathcal{M} \geq U$ coincides with (\ref{meas_gr_U}). To show that $\mathcal{M}$ majorizes the second term of~(\ref{M_ineq_omega}) we write the analogue of (\ref{delta_case}):
\begin{multline}
\label{delta_case_dom}
\int\limits_{\Omega} \varphi_{2\epsilon}^2 d\mathcal{M} 
=
\lim \limits_{k \to \infty} \int\limits_{0}^{+\infty} \int\limits_{\Omega} \Bigl[ t^{1-2s}|\nabla_{X}\left[w_k(X)\varphi_{2\epsilon}\right]|^2 
-
2 t^{1-2s}|\nabla_{x}w_k \cdot \nabla_{x}\varphi_{2\epsilon}\cdot w_k\varphi_{2\epsilon}|
\\+ 
t^{1-2s}|w_k(X) \nabla_{x}\varphi_{2\epsilon}(x)|^2 \Bigr] \, dX 
=: \widetilde{D}_1 - \widetilde{D}_2 + \widetilde{D}_3
\end{multline}
with $\widetilde{D}_1 \geq \mathcal{S}^{Sp}_{s, \sigma}(\Omega) \alpha_0^{\tfrac{2}{2^*_{\sigma}}}.$ The next step is the analogue of Lemma \ref{D3_lemma}. Indeed, we have
\begin{equation*}
\widetilde{D}_3 = \widetilde{a}_k +\widetilde{b}_k + \widetilde{c}_k := \left( \int\limits_{0}^{\delta} + \int\limits_{\delta}^{\frac{1}{\delta}} + \int\limits_{\frac{1}{\delta}}^{+\infty}\right) \int\limits_{\Omega} t^{1-2s} w^2_k(X) |\nabla_{x}\varphi_{2\epsilon}(x)|^2 \, dX.
\end{equation*}
We can pass to the limit in $\widetilde{b}_k,$ and the only remaining step is to obtain an analogue of~(\ref{a_k_c_k}). For a bounded $\Omega$ there is no explicit formula for the Green function, but we have the representation via Fourier series (see \cite[(3.1)-(3.8)]{Stinga}):
\begin{equation}
\label{fourier_rep_omega}
w(X) = \sum\limits_{i} d_{i}(t) \phi_i(x) \quad  \mbox{with} \quad d_{i}(t) = t^{s} \frac{2^{1-s}}{\Gamma(s)} \lambda^{s/2}_i \langle u, \phi_i \rangle \mathcal{K}_s(\lambda^{1/2}_i t),
\end{equation}
where $\mathcal{K}_s(\tau)$ is the modified Bessel function of the second kind; $\lambda_i, \phi_i$ were introduced in~(\ref{seminorm_Sp}). The asymptotic behavior of $\mathcal{K}_s$ is (see, e.g., \cite[(3.7)]{Stinga}):  
\begin{equation}
\label{bes_est_omega}
\mathcal{K}_s(\tau) \sim \Gamma(s) 2^{s-1} \tau^{-s} \ \quad \mbox{as} \quad \tau \to 0; \quad \mathcal{K}_s(\tau) \sim \left( \frac{\pi}{2\tau} \right)^{\frac{1}{2}} e^{-\tau} \left(1 + O(\tau^{-1}) \right) \quad \mbox{as} \quad \tau \to \infty.
\end{equation}
Thus, $w_k(X)$ can be estimated as (obviously, $\lambda_i \to \infty$)
\begin{equation*}
w_k(X) \leq C\sum\limits_{i} \langle u_k, \phi_i \rangle \phi_i(x) \cdot
\begin{cases}
1 &  \mbox{for} \quad  t \in [0, \delta]; \\
t^{2s-2} & \mbox{for} \quad t \in \left[\frac{1}{\delta}, +\infty\right),
\end{cases}
\end{equation*}
what gives 
\begin{equation*}
\widetilde{a}_k + \widetilde{c}_k 
\leq
\frac{C}{\epsilon^2} \int\limits_{0}^{\infty} \left( t^{1-2s} \chi_{[0, \delta]}(t) + t^{2s-3} \chi_{\left[\frac{1}{\delta}, +\infty\right)}(t) \right) \int\limits_{\Omega} u^2_k \ dX 
\leq
C(\epsilon) \delta^{2-2s} \|u_k\|_{L_2(\Omega)}^2.
\end{equation*}
Further, repeating argument from Sec. 5, we get (\ref{M_ineq_omega}). Similarly to (\ref{measure_ineq}) we have two alternatives: either $\alpha_0 = 0$ and the minimizer exists, or $\alpha_0 = 1$ and $u \equiv 0.$ We claim that in the second case the following inequality is fulfilled:
\begin{equation}
\label{second_case_est}
\mathcal{S}^{Sp}_{s, \sigma}(\Omega) \geq \mathcal{S}^{Sp}_{s, \sigma}(\mathbb{R}^n_+).
\end{equation}
Indeed, if $\{ u_k \}$ is a minimizing sequence for~(\ref{fu1}), then $\{ u_k \varphi_{2 \epsilon} \}$ is a minimizing sequence too: denominator of~(\ref{fu1}) converges to $\left[\alpha_0 \varphi_{2 \epsilon}^{2^*_{\sigma}}(0)\right]^{\frac{2}{2^*_{\sigma}}} = \alpha_0^{\frac{2}{2^*_{\sigma}}},$ while the convergence of numerator is controlled by~(\ref{delta_case_dom})
and Lemma \ref{D3_lemma} ($\widetilde{D}_2 = \widetilde{D}_3 = 0$):
\begin{equation*}
\lim \limits_{k \to \infty} \int\limits_{0}^{+\infty} \int\limits_{\Omega} t^{1-2s} \left|\nabla_{X}\left[w_k(X) \varphi_{2 \epsilon} \right]\right|^2 dX 
= 
\int\limits_{\Omega} \varphi_{2\epsilon}^2 d\mathcal{M}
=
\mathcal{S}^{Sp}_{s, \sigma}(\Omega)\alpha_0^{\frac{2}{2^*_{\sigma}}} \varphi_{2 \epsilon}^{2}(0) = \mathcal{S}^{Sp}_{s, \sigma}(\Omega)\alpha_0^{\frac{2}{2^*_{\sigma}}}.
\end{equation*}
Therefore, we can assume $u_k$ supported in $\mathbb{B}_{2 \epsilon}.$ Let $\Theta_{1}(x)$
be the coordinate transformation that flattens $\partial \Omega$ inside $\mathbb{B}_{r_0}:$
\begin{equation*}
y \equiv (y', y_n) = \Theta_{1}(x) := (x', x_n - F(x')) = x - F(x')e_n,
\end{equation*}
Jacobian of $\Theta_{1}(x)$ is equal to 1, thus
\begin{equation*}
\mathcal{I}_{\sigma, \Omega}[u_k] =  \frac{C_s \mathcal{E}_s \left[w_k\right]}{\| |x|^{\sigma-s} u_k  \|^2_{L_{2^*_{\sigma}}(\Omega)}} 
=
\frac{\int\limits_{0}^{+\infty} \int\limits_{\mathbb{R}^n_{+}} z^{1-2s}|\nabla_{Y} w_k(y', y_n + F(y'), z)|^2 \, dY \cdot (1 + o_{\epsilon}(1))}{\int\limits_{\mathbb{R}^n_{+}} \left||y'|^2 + (y_n + F(y'))^2\right|^{\frac{(\sigma - s)2^*_{\sigma}}{2}} \cdot u_k^{2^*_{\sigma}}(y', y_n + F(y')) \, dy}.
\end{equation*}
Since $w_k(y', y_n + F(y'), z)$ is an admissible extension of $u_k(y', y_n + F(y')),$ we have
\begin{equation*}
\mathcal{I}_{\sigma, \Omega}[u_k] \geq
\mathcal{S}^{Sp}_{s, \sigma}(\mathbb{R}^n_+) \cdot (1 + o_{\epsilon}(1)),
\end{equation*}
what gives (\ref{second_case_est}). 

To complete the proof we use the assumptions (\ref{f_convexity})-(\ref{f2_f_relations}) on $\partial \Omega$ to construct a function~$\Phi_{\varepsilon}(x)$ such that $\mathcal{I}_{\sigma, \Omega}[\Phi_{\varepsilon}(x)] < \mathcal{S}^{Sp}_{s, \sigma}(\mathbb{R}^n_+).$ We define $\Theta_{\varepsilon}(x)$ and $\Theta_{\varepsilon}(X)$ as 
\begin{equation}
\Theta_{\varepsilon}(x) := \varepsilon^{-1}  \Theta_{1}(x), \quad \Theta_{\varepsilon}(X) := \left(\Theta_{\varepsilon}(x),\varepsilon^{-1} t \right) = \left(\varepsilon^{-1}(x-F(x')e_n),\varepsilon^{-1} t \right). \label{theta_x}
\end{equation}
Jacobians of $\Theta_{\varepsilon}(x)$ and $\Theta_{\varepsilon}(X)$ are equal to $\varepsilon^{-n}$ and $\varepsilon^{-n-1}$ respectively. Let $\delta \in \left(0, r_0\right),$ we define  $\widetilde{\varphi}(x) :=\varphi_{\delta}\left(\Theta_{1}(x)\right).$ Note that $\widetilde{\varphi}(\Theta^{-1}_{\varepsilon}(y))$ is radial:
\begin{equation*}
\widetilde{\varphi}(\Theta^{-1}_{\varepsilon}(y)) = \widetilde{\varphi}(\varepsilon y', \varepsilon y_n + F(\varepsilon y')) =  \varphi_{\delta}(\Theta_{1}(\Theta^{-1}_{\varepsilon}(y))) = \varphi_{\delta}(\varepsilon|y|).
\end{equation*}
Now we put 
\begin{equation*}
\Phi_{\varepsilon}(x) := \varepsilon^{-\tfrac{n-2s}{2}} \Phi(\Theta_{\varepsilon}(x))\widetilde{\varphi}(x); \quad
w_{\varepsilon}(X) := \varepsilon^{-\tfrac{n-2s}{2}} \mathcal{W}(\Theta_{\varepsilon}(X))\widetilde{\varphi}(x)
\end{equation*}
(recall that $\Phi(y)$ is a minimizer of (\ref{fu1}) in $\mathbb{R}^n_+$ and $\mathcal{W}(Y)$ is its Stinga--Torrea extension). Obviously, $w_{\varepsilon}(X)$ is an admissible extension of $\Phi_{\varepsilon}(x),$ therefore
\begin{equation}
\label{eq3}
\mathcal{I}_{\sigma, \Omega}[\Phi_{\varepsilon}(x)]
= 
\frac{\langle (-\Delta)^s_{Sp} \Phi_{\varepsilon}, \Phi_{\varepsilon}\rangle}{\| |x|^{\sigma-s} \Phi_{\varepsilon}(x)  \|^2_{L_{2^*_{\sigma}}(\Omega)}} 
\leq
\frac{ \int \limits_0^{+\infty}  \int\limits_{\Omega} t^{1-2s} |\nabla_{X}  w_{\varepsilon}(X)|^2 \, dX}{\| |x|^{\sigma-s} \Phi_{\varepsilon}(x)  \|^2_{L_{2^*_{\sigma}}(\Omega)}}.
\end{equation}
In Secs. 8 and 9 we derive the estimates for the numerator and denominator in the right-hand side of (\ref{eq3}):
\begin{gather}
\label{denom_est}
\int\limits_{\Omega} \frac{|\Phi_{\varepsilon}(x)|^{2^*_{\sigma}}}{|x|^{(s-\sigma)2^*_{\sigma}}} \, dx = 1 - \mathcal{A}_1(\varepsilon)\cdot(1 + o_{\varepsilon}(1) + o_{\delta}(1)); \\
\label{num_est}
\mathcal{E}_s \left[ w_{\varepsilon} \right] = \mathcal{S}^{Sp}_{s, \sigma}(\mathbb{R}^n_+) + \mathcal{A}_2(\varepsilon) \cdot (1 + o_{\varepsilon}(1) + o_{\delta}(1)) - \frac{2\mathcal{S}^{Sp}_{s, \sigma}(\mathbb{R}^n_+)}{2^*_{\sigma}} \mathcal{A}_1(\varepsilon)\cdot(1 + o_{\varepsilon}(1)),
\end{gather}
where, for a fixed $\delta$ and $\varepsilon \to 0,$ 
\begin{equation*}
\mathcal{A}_1 (\varepsilon) \sim c_1 \varepsilon^{-1} f(\varepsilon); \quad
\mathcal{A}_2 (\varepsilon) \sim c_2 \varepsilon^{-1} f(\varepsilon); \quad
\mathcal{A}_1 (\varepsilon), \mathcal{A}_2 (\varepsilon) < 0
\end{equation*}
with $c_1, c_2 > 0.$ Therefore, for sufficiently small $\delta$ and $\varepsilon$ we have 
\begin{multline*}
\mathcal{I}_{\sigma, \Omega}[\Phi_{\varepsilon}(x)]
\leq
\frac{\mathcal{S}^{Sp}_{s, \sigma}(\mathbb{R}^n_+) + \mathcal{A}_2(\varepsilon) \cdot (1 + o_{\varepsilon}(1) + o_{\delta}(1)) - \frac{2\mathcal{S}^{Sp}_{s, \sigma}(\mathbb{R}^n_+)}{2^*_{\sigma}} \mathcal{A}_1(\varepsilon) \cdot (1 + o_{\varepsilon}(1))}{\Bigl(1 - \mathcal{A}_1(\varepsilon) \cdot (1 + o_{\varepsilon}(1) + o_{\delta}(1))\Bigr)^{\frac{2}{2^*_{\sigma}}}} 
\\=
\mathcal{S}^{Sp}_{s, \sigma}(\mathbb{R}^n_+) + \mathcal{A}_2(\varepsilon) \cdot (1 + o_{\varepsilon}(1) + o_{\delta}(1)) 
<
\mathcal{S}^{Sp}_{s, \sigma}(\mathbb{R}^n_+).
\end{multline*}
Thus (\ref{second_case_est}) is not fulfiled, the minimizer exists and the Theorem \ref{main_theorem} is proved.
\end{proof}
\section{Estimates for $\Phi(y)$ and $\mathcal{W}(Y)$}
This section is devoted to the proof of Lemma \ref{min_estimates}. As a first step, we obtain
the ``rough'' estimate for $\Phi(y)$ using the method from \cite[Lemma 3.5]{Scheglova} (see also \cite[Sec. II.5]{Lad_Ur}): it bounds~$\Phi$ in terms of its modulus of continuity in Lebesgue space with the critical Sobolev exponent. 

Let $0 < \tau < \|\Phi\|_{L_{2^*_{s}}\left(\mathbb{R}^n_{+}\right)},$ then there exists the \textbf{level} $\lambda := \lambda \left(\Phi, \tau\right)$ such that
\begin{equation*}
\|\Phi-\lambda\|_{L_{2^*_{s}} \left(\mathcal{Q}_{\lambda}\right)} = \tau, \quad \mbox{where} \quad\mathcal{Q}_{\lambda} := \{y \in \mathbb{R}^n_{+} : \Phi(y) > \lambda\}.
\end{equation*}
\begin{lemma}
\label{rough_phi}
There exists $\tau_*\left(n, s, \sigma\right)$ such that for any positive solution $\Phi(y)$ of (\ref{main_equation}) in $\mathbb{R}^n_{+}:$
\begin{equation}
\label{rough_zero}
\sup \Phi \leq C \cdot \lambda(\Phi, \tau_*).
\end{equation}
\end{lemma}
\begin{proof}
For any $\eta(Y) \in W^1_2(\mathbb{R}_+^n \times \mathbb{R}_+, z^{1-2s}),$ $\eta|_{y_n = 0} = 0$ we have
\begin{equation*}
\int\limits_{0}^{+\infty} \int\limits_{\mathbb{R}_+^n} z^{1-2s}  \nabla_{Y} \mathcal{W}(Y) \cdot \nabla_{Y} \eta(Y) \, dY  = \int\limits_{\mathbb{R}_{+}^n} \frac{\Phi^{2^*_{\sigma} - 1}(y)}{|y|^{(s-\sigma)2^*_{\sigma}}} \eta(y, 0) \, dy.
\end{equation*}
Let $\eta(Y) := \left[\mathcal{W}(Y) -  \lambda \right]_{+}$:
\begin{equation*}
\mathcal{E}_{s, \lambda}\left[\mathcal{W}\right] 
:= 
\int\limits_{\{\mathcal{W} > \lambda\}} z^{1-2s}  |\nabla_{Y} \mathcal{W}(Y)|^2 \, dY 
=
\int\limits_{\mathcal{Q}_{\lambda}}  \frac{\Phi^{2^*_{\sigma} - 1}(y)}{|y|^{(s-\sigma)2^*_{\sigma}}}  \left[\Phi(y) - \lambda \right] \, dy \leq \int\limits_{\mathcal{Q}_{\lambda}}  \frac{\Phi^{2^*_{\sigma}}(y)}{|y|^{(s-\sigma)2^*_{\sigma}}} \, dy.
\end{equation*}
Let us estimate the integral from the right-hand side:
\begin{multline*}
\|\Phi\|^{2^*_{\sigma}}_{L_{2^*_{\sigma}} \left(\mathcal{Q}_{\lambda}, |y|^{(\sigma-s)2^*_{\sigma}}\right)} \leq \left( \|\Phi - \lambda\|_{L_{2^*_{\sigma}} \left(\mathcal{Q}_{\lambda}, |y|^{(\sigma-s)2^*_{\sigma}}\right)} + \lambda \|1\|_{L_{2^*_{\sigma}} \left(\mathcal{Q}_{\lambda}, |y|^{(\sigma-s)2^*_{\sigma}} \right)}  \right)^{2^*_{\sigma}} 
\\ \leq 
2^{2^*_{\sigma}} \left( \|\Phi - \lambda\|^{2^*_{\sigma}}_{L_{2^*_{\sigma}}} + \lambda ^{2^*_{\sigma}} \int \limits_{\mathcal{Q}_{\lambda}} \frac{1}{|y|^{(s- \sigma)2^*_{\sigma}}} \, dy \right)
\stackrel{*}{\leq}
2^{2^*_{\sigma}} \|\Phi - \lambda\|^{2^*_{\sigma}}_{L_{2^*_{\sigma}} \left(\mathcal{Q}_{\lambda}, |y|^{(\sigma-s)2^*_{\sigma}}\right)} + C_1 {\lambda}^{2^*_{\sigma}} |\mathcal{Q}_{\lambda}|^{\frac{n-2s}{n-2\sigma}},
\end{multline*}
the inequality (*) follows from the Schwarz symmetrization. Recall that $\tau \equiv \|\Phi - \lambda\|_{L_{2^*_{s}}\left(\mathcal{Q}_{\lambda}\right)},$ using the H{\"o}lder inequality we get
\begin{multline*}
\|\Phi - \lambda\|^{2^*_{\sigma}}_{L_{2^*_{\sigma}} \left(\mathcal{Q}_{\lambda}, |y|^{(\sigma-s)2^*_{\sigma}}\right)} 
\leq 
\|\Phi - \lambda\|^{\frac{2n(s-\sigma)}{(n-2\sigma)s}}_{L_2 \left(\mathcal{Q}_{\lambda}, |y|^{-2s} \right)} \cdot \|\Phi - \lambda\|^{2 - \frac{2n(s-\sigma)}{(n-2\sigma)s} + 2^*_{\sigma} - 2}_{L_{2^*_{s}}\left(\mathcal{Q}_{\lambda}\right)}
\\=
\|\Phi - \lambda\|^{\frac{2n(s-\sigma)}{(n-2\sigma)s}}_{L_2 \left(\mathcal{Q}_{\lambda}, |y|^{-2s} \right)} \cdot \|\Phi - \lambda\|^{2 - \frac{2n(s-\sigma)}{(n-2\sigma)s}}_{L_{2^*_{s}}\left(\mathcal{Q}_{\lambda}\right)} \cdot \tau^{2^*_{\sigma} - 2}.
\end{multline*}
Due to fractional Hardy and Sobolev inequaities
\begin{equation*}
\|\Phi - \lambda\|^{\frac{2n(s-\sigma)}{(n-2\sigma)s}}_{L_2 \left(\mathcal{Q}_{\lambda}, |y|^{-2s} \right)} \cdot \|\Phi - \lambda\|^{2 - \frac{2n(s-\sigma)}{(n-2\sigma)s}}_{L_{2^*_{s}}\left(\mathcal{Q}_{\lambda}\right)}
\leq
C_2 \langle (-\Delta)^s_{\mathcal{Q}_{\lambda}, Sp} [\Phi - \lambda]_+, [\Phi - \lambda]_+ \rangle
\stackrel{**}{\leq}
C_2 \mathcal{E}_{s, \lambda}\left[\mathcal{W}\right],
\end{equation*}
the inequality (**) follows from the fact that $\eta(Y)$ is an admissible extension of $[\Phi - \lambda]_+.$ To sum up,
\begin{equation*}
\mathcal{E}_{s, \lambda}\left[\mathcal{W}\right] 
\leq
{2^*_{\sigma}} C_2 \mathcal{E}_{s, \lambda}\left[\mathcal{W}\right] \tau^{2^*_{\sigma} - 2}  + C_1 {\lambda}^{2^*_{\sigma}} |\mathcal{Q}_{\lambda}|^{\frac{n-2s}{n-2\sigma}}.
\end{equation*}
Suppose that $\tau_{*}$ satisfies $2^{2^*_{\sigma}}C_2 \tau_{*}^{2^*_{\sigma} - 2} \leq \frac{1}{2}.$ For all $\lambda > \lambda(\Phi, \tau_{*})$ we have
\begin{equation}
\label{ineq_meas}
C_3 \|\Phi - \lambda\|^{2}_{L_{2^*_{s}} \left(\mathcal{Q}_{\lambda}\right)} \leq \mathcal{E}_{s, \lambda}\left[\mathcal{W}\right] \leq 2C_1 {\lambda}^{2^*_{\sigma}} |\mathcal{Q}_{\lambda}|^{\frac{n-2s}{n-2\sigma}}.
\end{equation}
From (\ref{ineq_meas}) we obtain
\begin{equation}
\label{lambda_ineq}
\mathfrak{g}(\lambda) 
:= 
\int\limits_{\mathcal{Q}_{\lambda}} \left[\Phi(y)-\lambda\right] \, dy 
\leq
\|\Phi - \lambda\|_{L_{2^*_{s}} \left(\mathcal{Q}_{\lambda}\right)} \cdot |\mathcal{Q}_{\lambda}|^{\frac{n+2s}{2n}}
\leq
C_4 {\lambda}^{\frac{n}{n-2\sigma}} |\mathcal{Q}_{\lambda}|^{1 + \frac{\sigma (n-2s)}{n(n-2\sigma)}}.
\end{equation}
Using the layer cake representation for the Lebesgue integral
\begin{equation*}
\mathfrak{g}(\lambda)
=
\int\limits_{\mathcal{Q}_{\lambda}} \int\limits_{\lambda}^{\infty} \chi_{\{\theta < \Phi(y)\}} \, d\theta \, dy
= 
\int\limits_{\lambda}^{\infty} |\mathcal{Q}_{\theta}| \, d\theta
\end{equation*}
we get $\mathfrak{g}'(\lambda) = - |\mathcal{Q}_{\lambda}|$ for a.e. $\lambda.$ Thus (\ref{lambda_ineq}) takes the form
\begin{equation*}
-\mathfrak{g}'(\lambda)\left[\mathfrak{g}(\lambda)\right]^{- \frac{n(n-2\sigma)}{n^2 - n\sigma - 2 \sigma s}}
\geq
C_5 \lambda^{- \frac{n^2}{n^2 - n\sigma - 2 \sigma s}}.
\end{equation*}
By integrating over the segment $[\lambda, \sup\Phi]$ we get
\begin{gather*}
-\mathfrak{g}(\lambda)^{\frac{n\sigma - 2 \sigma s}{n^2 - n\sigma - 2 \sigma s}}
\leq
C_6 \left[(\sup \Phi)^{- \frac{n\sigma + 2 \sigma s}{n^2 - n\sigma - 2 \sigma s}} - \lambda^{- \frac{n\sigma + 2 \sigma s}{n^2 - n\sigma - 2 \sigma s}} \right]; \\
(\sup \Phi)^{- \frac{n\sigma + 2 \sigma s}{n^2 - n\sigma - 2 \sigma s}}
\geq 
\lambda^{- \frac{n\sigma + 2 \sigma s}{n^2 - n\sigma - 2 \sigma s}} - C^{-1}_6 \mathfrak{g}(\lambda)^{\frac{n\sigma - 2 \sigma s}{n^2 - n\sigma - 2 \sigma s}}.
\end{gather*}
Using (\ref{lambda_ineq}) for $\tau_{*} \leq \left( \frac{C_6}{2} \right)^{\frac{n^2 - n\sigma - 2 \sigma s}{n\sigma - 2 \sigma s}} \cdot \left[\mathcal{S}_{s,s}^{-1} \cdot \mathcal{S}^{Sp}_{s, \sigma}(\mathbb{R}^n_{+})\right]^{-\frac{n+2s}{2(n-2s)}}$ we obtain
\begin{multline}
\label{g_ineq}
\mathfrak{g}(\lambda)\lambda^{\frac{n+2s}{n-2s}} 
\leq\|\Phi - \lambda\|_{L_{2^*_{s}} \left(\mathcal{Q}_{\lambda} \right)} \left( |\mathcal{Q}_{\lambda}| \lambda^{2^*_s} \right)^{\frac{n+2s}{2n}}
\leq 
\tau_{*} \|\Phi\|^{2^*_{s}-1}_{L_{2^*_{s}}}
\\ \leq
\tau_{*}  \left[ \mathcal{S}_{s,s}^{-1} \mathcal{E}_{s}\left[\mathcal{W}\right]\right]^{\frac{n+2s}{2(n-2s)}}
=
\tau_{*}  \left[ \mathcal{S}_{s,s}^{-1} \cdot \mathcal{S}^{Sp}_{s, \sigma}(\mathbb{R}^n_{+})\right]^{\frac{n+2s}{2(n-2s)}}
\leq 
\left( \frac{C_6}{2} \right)^{\frac{n^2 - n\sigma - 2 \sigma s}{n\sigma - 2 \sigma s}},
\end{multline}
what gives the required inequality (\ref{rough_zero}):
\begin{equation*}
(\sup \Phi)^{- \frac{n\sigma + 2 \sigma s}{n^2 - n\sigma - 2 \sigma s}}
\geq 
\frac{1}{2} \lambda^{- \frac{n\sigma + 2 \sigma s}{n^2 - n\sigma - 2 \sigma s}}. \qedhere
\end{equation*}
\end{proof}
\begin{corollary}
Any minimizer $\Phi(y)$ admits the estimate ($\tau_{*}$ was introduced in Lemma \ref{rough_phi}):
\begin{equation}
\label{rough_infty}
\Phi(y) \leq \frac{C(n,s, \sigma, \lambda\left(\Phi, \tau_{*}\right), \lambda\left(\Phi^*, \tau_{*}\right))}{(1 + |y|)^{n-2s}}.
\end{equation} 
\end{corollary}
\begin{proof}
For $|y| \leq 1$ the estimate (\ref{rough_infty}) coincides with (\ref{rough_zero}). For $|y| \geq 1$ the estimate (\ref{rough_infty}) can be obtained via the $s$-Kelvin transform (\ref{s_kelvin}):
\begin{equation*}
\Phi (y) \leq \frac{1}{|y|^{n-2s}} \cdot \sup \Phi^*\left( \frac{y}{|y|^2} \right) \leq \frac{C(n,s, \sigma, \lambda\left(\Phi^*,  \tau_{*}\right))}{|y|^{n-2s}}. \qedhere
\end{equation*}
\end{proof}
\begin{proof}[\textbf{Proof of Lemma \ref{min_estimates}}]
The estimate for $\Phi(y)$ in  (\ref{W_common_ineq}) follows from the estimate for $\mathcal{W}(Y)$ due to $\Phi(y) = \mathcal{W}(y,0).$ Moreover,
$s$-Kelvin transform argument shows that it suffices to prove~(\ref{W_common_ineq}) for $|Y| \leq 1$ only. Using the Green function (\ref{green_w_h_Rn+}), we can write
\begin{equation*}
\mathcal{W}(Y) 
=
\Biggl( \ \int\limits_{|\xi| > 2}
+
\int\limits_{\substack{|\xi| \leq 2 \\ |y - \xi| > \frac{y_n}{2}}}
+
\int\limits_{\substack{|\xi| \leq 2 \\ |y - \xi| \leq \frac{y_n}{2}}} \Biggr) G_s(Y, \xi) \, \frac{\Phi^{2^*_{\sigma} - 1}(\xi)}{|\xi|^{(s-\sigma)2^*_{\sigma}}} \, d\xi
=:
A_1 + A_2 + A_3.
\end{equation*}
To estimate $A_1$, we use (\ref{rough_infty}) and (\ref{green_est}) with $\mathfrak{b} = 1:$
\begin{equation}
\label{I1_ineq}
A_1 
\leq 
C y_n \int\limits_{|\xi| > 2} |\xi|^{(2^*_{\sigma} - 1)(2s - n)} |\xi|^{(\sigma - s)2^*_{\sigma}} \frac{\xi_n}{|\xi|^{n-2s+2}} \, d\xi 
\leq
C y_n \int\limits_{|\xi| > 2} |\xi|^{-\left(\frac{n^2  - 2ns}{n - 2\sigma} + n + 1\right)} \, d\xi \leq C y_n.
\end{equation}
Estimates of $A_2$ and $A_3$ are obtained iteratively. Recall that we have fixed the minimizer~$\Phi(y).$ Let the following a priori estimate with $\mathfrak{p} \in [0, 1)$ be fulfilled (for $\mathfrak{p} = 0$ it was proved in Lemma \ref{rough_phi}):
\begin{equation}
\label{Phi_recursive}
\Phi(y) \leq C y_n^{\mathfrak{p}}.
\end{equation}
We claim that (\ref{Phi_recursive}) implies 
\begin{equation}
\label{Phi_recursive_p*}
\mathcal{W}(Y) \leq C y_n^{\mathfrak{p}_*} \quad \mbox{and} \quad \Phi(y) \leq C y_n^{\mathfrak{p}_*},
\end{equation}
with $\mathfrak{p}_* := \min(\mathfrak{q} + \mathfrak{p}, 1)$ and
\begin{equation*}
\mathfrak{q} :=\frac{\sigma(n-2s)}{n-2\sigma} =  s - \frac{(s - \sigma) 2^*_{\sigma}}{2} \in (0, s).
\end{equation*}
Indeed, to estimate $A_2$ we notice that on the integration set one has
\begin{equation*}
\xi_n \leq |\xi-y| + y_n \leq 3|\xi-y|,
\end{equation*}
therefore inequalities $\xi_n < |\xi|,$ (\ref{rough_infty}) and (\ref{green_est}) with $\mathfrak{b} = 1$ give us
\begin{equation}
\label{I2_ineq}
A_2 
\leq
\int\limits_{\substack{|\xi| \leq 2 \\ |y - \xi| > \frac{y_n}{2}}} |\xi|^{(\sigma-s)2^*_{\sigma}} \frac{C y_n\xi_n^{1+(2^*_{\sigma} - 1)\mathfrak{p}} }{|y - \xi|^{n-2s+2}} \, d\xi 
\leq
Cy_n^{\mathfrak{p}_*} \int\limits_{|\xi| \leq 2} |\xi|^{2(\mathfrak{q} - s)}
 |y - \xi|^{-n + 2s + (2^*_{\sigma} - 1)\mathfrak{p} - \mathfrak{p}_*} \, d\xi,
\end{equation}
both exponents are negative, their sum is greater than $-n$ and the integral converges.

To estimate $A_3,$ we notice that on the integration set one has
\begin{equation*}
|\xi| \geq |y| - |y - \xi| \geq |y| - \frac{y_n}{2} \geq \frac{y_n}{2}; \quad \xi_n \leq |y_n - \xi_n| + y_n \leq \frac{3 y_n}{2},
\end{equation*}
therefore (\ref{green_est}) with $\mathfrak{b} = 0$ gives us
\begin{equation}
\label{I3_ineq}
A_3
\leq
\int\limits_{|y - \xi| \leq \frac{y_n}{2}} \frac{C \xi^{(2^*_{\sigma} - 1)\mathfrak{p}}_n}{|\xi|^{2s - 2\mathfrak{q}}|y - \xi|^{n-2s}} \, d\xi 
\leq
C y_n^{(2^*_{\sigma} - 1)\mathfrak{p} - 2s + 2\mathfrak{q}} \int\limits_{|y - \xi| \leq \frac{y_n}{2}} \frac{1}{|y - \xi|^{n-2s}} \, d\xi 
\leq
C y_n^{\mathfrak{p}_*}.
\end{equation}
Putting (\ref{I1_ineq}), (\ref{I2_ineq}) and (\ref{I3_ineq}) together, we obtain (\ref{Phi_recursive_p*}), i.e. we have increased the exponent in~(\ref{Phi_recursive}) by at least   $\min(\mathfrak{q}, 1-\mathfrak{p}).$ Iterating this process, we get (\ref{Phi_recursive_p*}) with $\mathfrak{p}_* = 1.$ The estimate (\ref{W_common_ineq}) is proved completely.

To prove (\ref{GradW_common_ineq}) we have to derive estimates at the origin and at infinity separately because $\mathcal{V}(y)$ is not invariant under the $s$-Kelvin transform. For $|y| \leq 1,$ we write the integral representation for $\nabla_{Y} \mathcal{W}(Y)$ as follows
\begin{equation*}
\nabla_{Y} \mathcal{W}(Y)
=
\left( \ \int\limits_{|\xi| \geq 2} +  \int\limits_{|\xi| < 2} \right) \frac{\Phi^{2^*_{\sigma} - 1}(\xi)}{|\xi|^{(s-\sigma)2^*_{\sigma}}} \nabla_{Y} G_s(Y, \xi) \, d\xi
=:
A_4 + A_5.
\end{equation*}
Obviously,
\begin{equation}
\label{V_estimate}
\mathcal{V}(y) 
\leq
2 \int\limits_{0}^{+\infty} z^{1-2s} A^2_4(Y) \, dz +2 \int\limits_{0}^{2} z^{1-2s} A^2_5(Y) \, dz +2 \int\limits_{2}^{+\infty} z^{1-2s} A^2_5(Y) \, dz.
\end{equation}
We estimate $A_4$ using (\ref{green_grad_est}) and (\ref{W_common_ineq}):
\begin{equation*}
A_4 \leq  \int\limits_{|\xi| \geq 2} \frac{\xi^{2^*_{\sigma} - 1}_n}{|\xi|^{(s-\sigma)2^*_{\sigma}+ (2^*_{\sigma}-1)(n-2s+2)} \left(|y - \xi|^2 + z^2\right)^{\frac{n-2s+1}{2}}} \, d\xi.
\end{equation*}
Therefore, taking into account $|y - \xi| \geq |\xi| - |y| \geq 1$ we get
\begin{equation*}
\int\limits_{0}^{+\infty} z^{1-2s} A^2_4(Y) \, dz
\leq
C \int\limits_{0}^{+\infty} \frac{z^{1-2s}}{(1 + z^2)^{n-2s+1}} \, dz \cdot \left( \int\limits^{+\infty}_{2} \frac{r^{2^*_{\sigma} - 1} r^{n - 1}}{r^{(s-\sigma)2^*_{\sigma} + (2^*_{\sigma}-1)(n-2s+2)}} \, dr \right)^2 \leq C,
\end{equation*}
convergence of the last integral follows from the equality
\begin{equation*}
2^*_{\sigma} - 1 + n - 1 - (s-\sigma)2^*_{\sigma} - (2^*_{\sigma}-1)(n-2s+2) = -\frac{2\sigma(n-2s+2)}{n-2\sigma} - 2.
\end{equation*}
The estimate of $A_5$ also follows from (\ref{green_grad_est}) and (\ref{W_common_ineq}):
\begin{equation*}
A_5 \leq  \int\limits_{|\xi| < 2} \frac{\xi^{2^*_{\sigma} - 1}_n}{|\xi|^{(s-\sigma)2^*_{\sigma}} \left(|y - \xi|^2 + z^2\right)^{\frac{n-2s+1}{2}}} \, d\xi.
\end{equation*}
Using this inequality, we estimate the second term in (\ref{V_estimate}):
\begin{equation*}
\int\limits_{0}^{2} z^{1-2s} A^2_5(Y) \, dz
\leq
 C \int\limits_{0}^{2} z^{-1+\min(s, 1-s)} \, dz \cdot \left( \ \int\limits_{|\xi| < 2} \frac{\xi^{2^*_{\sigma} - 1}_n }{|\xi|^{(s-\sigma)2^*_{\sigma}} |y - \xi|^{n-s+\frac{\min(s, 1-s)}{2}}} \, d\xi \right)^2 \leq C,
\end{equation*}
convergence of the last integral follows from the inequality
\begin{equation*}
2^*_{\sigma} - 2 - (s-\sigma)2^*_{\sigma} + s - \frac{\min(s, 1-s)}{2} = \frac{2\sigma(n - 2s + 2)}{n-2\sigma} - s -  \frac{\min(s, 1-s)}{2} > -1.
\end{equation*}
Finally, the third term in (\ref{V_estimate}) can be estimated as
\begin{equation*}
\int\limits_{2}^{+\infty} z^{1-2s} A^2_5(Y) \, dz
\leq
 C \int\limits_{2}^{+\infty} \frac{z^{1-2s}}{z^{2n-4s+2}} \, dz \cdot \left( \int\limits^{2}_{0} \frac{r^{2^*_{\sigma} - 1} r^{n - 1}}{r^{(s-\sigma)2^*_{\sigma}}} \, dr \right)^2 \leq C,
\end{equation*}
and (\ref{GradW_common_ineq}) is proved for $|y| \leq 1.$

For $|y| > 1,$ we write the integral representation for $\nabla_{Y} \mathcal{W}(Y)$ as follows
\begin{equation*}
\nabla_{Y} \mathcal{W}(Y)
=
\left( \int\limits_{|y-\xi| < \frac{|y|}{10}} +  \int\limits_{|y-\xi| \geq \frac{|y|}{10}} \right) \frac{\Phi^{2^*_{\sigma} - 1}(\xi)}{|\xi|^{(s-\sigma)2^*_{\sigma}}} \nabla_{Y} G_s(Y, \xi) \, d\xi
=:
A_6 + A_7.
\end{equation*}
Then $\mathcal{V}(y)$ can be estimated with an obvious inequality
\begin{equation}
\label{V_estimate2}
\mathcal{V}(y)
\leq 
2\int\limits_{0}^{+\infty} z^{1-2s} A^2_6(Y) \, dz + 2\int\limits_{0}^{+\infty} z^{1-2s} A^2_7(Y) \, dz.
\end{equation}
We estimate $A_6$ using (\ref{green_grad_est}), (\ref{W_common_ineq}) and $|\xi| \geq \frac{9|y|}{10} \geq \frac{9}{10}:$ 
\begin{equation*}
A_6
\leq
C \int\limits_{|y-\xi| < \frac{|y|}{10}} \frac{|\xi|^{(\sigma-s)2^*_{\sigma} - (2^*_{\sigma}-1)(n-2s+1)}}{\left(|y - \xi|^2 + z^2\right)^{\frac{n-2s+1}{2}}} \, d\xi
\leq
\frac{C}{|y|^{(s-\sigma)2^*_{\sigma} + (2^*_{\sigma}-1)(n-2s+1)}} \int\limits_{0}^{\frac{|y|}{10}} \frac{r^{n-1}} {\left(r^2 + z^2\right)^{\frac{n-2s+1}{2}}} \, dr.
\end{equation*}
Changing the variable shows that
\begin{equation*}
\int\limits_{0}^{+\infty} z^{1-2s} \left(\int\limits_{0}^{\frac{|y|}{10}} \frac{r^{n-1}} {\left(r^2 + z^2\right)^{\frac{n-2s+1}{2}}} \, dr\right)^2 dz 
= 
\frac{|y|^{2s}}{10^{2s}} \int\limits_{0}^{+\infty} z^{1-2s} \left( \int\limits_{0}^{1} \frac{r^{n-1}} {\left(r^2 + z^2\right)^{\frac{n-2s+1}{2}}} \, dr \right)^2 dz,
\end{equation*}
what gives the estimate of the first term in (\ref{V_estimate2})
\begin{equation*}
\int\limits_{0}^{+\infty} z^{1-2s} A^2_6(Y) \, dz 
\leq
\frac{C|y|^{2s}}{|y|^{2((s-\sigma)2^*_{\sigma} + (2^*_{\sigma}-1)(n-2s+1))}}
=
\frac{C}{|y|^{2n - 2s + 2 + 4\sigma\frac{n-2s+2}{n-2\sigma}}}.
\end{equation*}
Finally, we estimate $A_7$ using (\ref{green_grad_est}) and (\ref{W_common_ineq}):
\begin{multline*}
A_7
\leq
\int\limits_{|y-\xi| \geq \frac{|y|}{10}} \frac{C \Phi^{2^*_{\sigma} - 1}(\xi)\xi_n}{|\xi|^{(s-\sigma)2^*_{\sigma}} \left(|y - \xi|^2 + z^2\right)^{\frac{n-2s+2}{2}}} \left( \frac{y_n}{\sqrt{|y - \xi|^2 + z^2}} + 1 \right) \, d\xi
\\ \leq
\frac{C}{\left(|y|^2 + z^2\right)^{\frac{n-2s+2}{2}}} \cdot \int\limits_{|y-\xi| \geq \frac{|y|}{10}} \frac{\xi^{2^*_{\sigma}}_n} {|\xi|^{(s-\sigma)2^*_{\sigma}} (1 + |\xi|^{(2^*_{\sigma}-1)(n-2s+2)}) } \, d\xi.
\end{multline*}
convergence of the last integral follows from the inequality
\begin{equation*}
2^*_{\sigma} - (s-\sigma)2^*_{\sigma} - (2^*_{\sigma}-1)(n-2s+2) =
- \frac{n^2  - 4s\sigma + 4\sigma}{n-2\sigma} = -n - 2\sigma \frac{n - 2s + 2}{n-2\sigma} < -n.
\end{equation*}
This gives the estimate of the second term in (\ref{V_estimate2})
\begin{equation*}
\int\limits_{0}^{+\infty} z^{1-2s} A^2_7(Y) \, dz
\leq 
C \int\limits_{0}^{+\infty} \frac{z^{1-2s} }{\left(|y|^2 + z^2\right)^{n-2s+2}} \, dz
\leq
\frac{C}{|y|^{2n-2s+2}},
\end{equation*}
and the estimate (\ref{GradW_common_ineq}) is proved completely!
\end{proof}
\section{Estimate of the denominator and derivation of (\ref{denom_est})}
To get (\ref{denom_est}) we modify the calculations from \cite[Sec. 4]{Dem_Naz}. We use the change of variables (\ref{theta_x}) and obtain by the Taylor formula:
\begin{multline*}
\int\limits_{\Omega} \frac{|\Phi_{\varepsilon}(x)|^{2^*_{\sigma}}}{|x|^{(s-\sigma)2^*_{\sigma}}} \, dx 
= 
\int\limits_{\mathbb{R}^n_+} \frac{|\Phi(y)|^{2^*_{\sigma}}}{|y + \varepsilon^{-1} F(\varepsilon y') e_n |^{(s-\sigma)2^*_{\sigma}}} \widetilde{\varphi}^{2^*_{\sigma}}(\Theta^{-1}_{\varepsilon}(y)) \, dy 
\\=
\int\limits_{\mathbb{R}^n_+} \frac{|\Phi(y)|^{2^*_{\sigma}}}{|y|^{(s-\sigma)2^*_{\sigma}}} \varphi_{\delta}^{2^*_{\sigma}}(\varepsilon y) \cdot 
\left( 1 - \frac{(s-\sigma)2^*_{\sigma}}{\varepsilon} F(\varepsilon y') \frac{y_n}{|y|^2} + \frac{F^2(\varepsilon y')}{\varepsilon^2 |y|^2} \cdot O_{\delta}(1) \right) dy 
\\= 
\int\limits_{\mathbb{R}^n_+} \frac{|\Phi(y)|^{2^*_{\sigma}}}{|y|^{(s-\sigma)2^*_{\sigma}}} \, dy 
-
\int\limits_{\mathbb{R}^n_+} \frac{|\Phi(y)|^{2^*_{\sigma}}}{|y|^{(s-\sigma)2^*_{\sigma}}} \left(1 - \varphi_{\delta}^{2^*_{\sigma}}(\varepsilon y) \right) dy
-
\int\limits_{\mathbb{R}^n_+} \frac{(s-\sigma)2^*_{\sigma} |\Phi(y)|^{2^*_{\sigma}}\varphi_{\delta}^{2^*_{\sigma}}(\varepsilon y) y_n}{\varepsilon |y|^{(s-\sigma)2^*_{\sigma} + 2}} F(\varepsilon y') \, dy 
\\+
O_{\delta}(1)
\int\limits_{\mathbb{R}^n_+} \frac{|\Phi(y)|^{2^*_{\sigma}}}{|y|^{(s-\sigma)2^*_{\sigma}}} \varphi_{\delta}^{2^*_{\sigma}}(\varepsilon y)  \frac{F^2(\varepsilon y')}{\varepsilon^2 |y|^2} \, dy 
=: I_1 - I_2- I_3 + I_4.
\end{multline*}
\begin{lemma}
\label{lemma_Ik}
The following relations hold:
\begin{enumerate}
\item $I_1 = 1$ and $I_2 \leq C \left( \frac{\varepsilon}{\delta} \right)^{\frac{n(n-2s+2)}{n-2\sigma}};$ 
\item 
\begin{equation}
\label{I3_estim}
\lim\limits_{\varepsilon \to 0} \varepsilon \frac{I_3}{f(\varepsilon)} = C \int\limits_{0}^{+\infty} \tau^{\alpha + n} \int\limits_{0}^{+\infty} \frac{|\Phi(\tau, \tau \varsigma)|^{2^*_{\sigma}} \varsigma d\varsigma}{|\tau^2 + \tau^2 \varsigma^2|^{\frac{(s-\sigma)2^*_{\sigma} + 2}{2}}} \, d\tau < +\infty;
\end{equation}
\item $\lim\limits_{\varepsilon \to 0} \left|\varepsilon \frac{I_4}{f(\varepsilon)} \right|  = o_{\delta}(1).$
\end{enumerate}
\end{lemma}
\begin{proof}
1. The equality $I_1 = 1$ is just a normalizing condition for $\Phi(y).$ Further, (\ref{W_common_ineq}) gives
\begin{equation*}
I_2 
\equiv 
\int\limits_{\mathbb{R}^n_+} \frac{|\Phi(y)|^{2^*_{\sigma}}}{|y|^{(s-\sigma)2^*_{\sigma}}} \left(1 - \varphi_{\delta}^{2^*_{\sigma}}(\varepsilon y) \right) dy
\leq
C \int \limits^{+\infty}_{\frac{\delta}{2\varepsilon}} r^{n-1-2^*_{\sigma}(n-s - \sigma +1)} \, dr
=
C\left(\frac{\varepsilon}{\delta}\right)^{\frac{n(n - 2s + 2)}{n-2\sigma}}.
\end{equation*}
2. We note that
\begin{multline*}
\frac{\varepsilon I_3}{f(\varepsilon)} 
=
\frac{(s-\sigma)2^*_{\sigma}}{f(\varepsilon)} \int\limits_{\mathbb{R}^n_+} \varphi_{\delta}^{2^*_{\sigma}}(\varepsilon y) \frac{|\Phi(y)|^{2^*_{\sigma}} y_n}{|y|^{(s-\sigma)2^*_{\sigma} + 2}} F(\varepsilon y') \, dy 
\\=
\frac{C}{f(\varepsilon)} \int\limits_{0}^{+\infty} \int\limits_{0}^{+\infty} \varphi_{\delta}^{2^*_{\sigma}}(\varepsilon \sqrt{\tau^2 + y_n^2}) \frac{|\Phi(\tau, y_n)|^{2^*_{\sigma}} y_n}{(\tau^2 + y_n^2)^{\frac{(s-\sigma)2^*_{\sigma} + 2}{2}}} \int\limits_{\mathbb{S}_{\tau}^{n-2}} F(\varepsilon y') \, d\mathbb{S}_{\tau}^{n-2}(y') \, dy_n \, d\tau
\\= 
C \int\limits_{0}^{+\infty} \tau^n \frac{f(\varepsilon \tau)}{f(\varepsilon)} \int\limits_{0}^{+\infty} \varphi_{\delta}^{2^*_{\sigma}}(\varepsilon \sqrt{\tau^2 + \tau^2\varsigma^2}) \frac{|\Phi(\tau, \tau \varsigma)|^{2^*_{\sigma}} \varsigma \, d\varsigma}{(\tau^2 + \tau^2\varsigma^2)^{\frac{(s-\sigma)2^*_{\sigma} + 2}{2}}} \, d\tau
\\= 
C\int\limits_{0}^{+\infty} \tau^{\alpha + n} \frac{\psi(\varepsilon \tau)}{\psi(\varepsilon)} \int\limits_{0}^{+\infty} \varphi_{\delta}^{2^*_{\sigma}}(\varepsilon \sqrt{\tau^2 + \tau^2\varsigma^2}) \frac{|\Phi(\tau, \tau\varsigma)|^{2^*_{\sigma}} \varsigma \, d\varsigma}{(\tau^2 + \tau^2\varsigma^2)^{\frac{(s-\sigma)2^*_{\sigma} + 2}{2}}} \, d\tau
=:
C\int\limits_{0}^{+\infty} P_{\varepsilon}(\tau) \, d\tau.
\end{multline*}
The pointwise limit of $P_{\varepsilon}(\tau)$ as $\varepsilon \to 0$ coincides with the integrand in the right-hand side of~(\ref{I3_estim}). To get the final result we use the Lebesgue dominated convergence theorem. To construct a summable majorant for $P_{\varepsilon}(\tau)$ we notice that $\psi(\tau)$ is an SVF and therefore $\psi(\tau) \tau^ {\beta}$ increases and $\psi(\tau) \tau^ {-\beta}$ decreases in the neighbourhood of the origin for any $\beta > 0,$ see \cite[Sec 1.5, (1)-(2)]{Seneta}. This implies
\begin{multline}
\label{psi_equal}
\chi_{\left[0, \frac{\delta}{\varepsilon}\right]}(\tau)
\frac{\psi(\varepsilon\tau)}{\psi(\varepsilon)} 
=
\frac{\psi(\varepsilon\tau) (\varepsilon\tau)^{\beta}}{\psi(\varepsilon) (\varepsilon)^{\beta}}
\chi_{[0,1]}(\tau) \tau^{-\beta} +\frac{\psi(\varepsilon\tau) (\varepsilon\tau)^{-\beta}}{\psi(\varepsilon) (\varepsilon)^{-\beta}}
\chi_{\left[1, \frac{\delta}{\varepsilon}\right]}(\tau) \tau^{\beta}
\\ \leq
C(\delta) \left(\chi_{[0,1]}(\tau) \tau^{-\beta} + \chi_{[1, +\infty)}(\tau) \tau^{\beta}\right).
\end{multline}
Thus,
\begin{equation*}
P_{\varepsilon}(\tau)
\leq
C(\delta) \left(\chi_{[0,1]}(\tau) \tau^{\alpha+n-\beta} + \chi_{[1, +\infty)}(\tau) \tau^{\alpha+n+\beta}\right)
\int\limits_{0}^{+\infty} \frac{|\Phi(\tau, \tau\varsigma)|^{2^*_{\sigma}} \varsigma}{(\tau^2 + \tau^2\varsigma^2)^{\frac{(s-\sigma)2^*_{\sigma} + 2}{2}}} \, d\varsigma.
\end{equation*}
By (\ref{W_common_ineq}), for $\tau \in  [0,1]$ we have
\begin{equation*}
\int\limits_{0}^{+\infty} \frac{(\tau^2 + \tau^2\varsigma^2)^{ \frac{(1 - s + \sigma)2^*_{\sigma} - 2}{2}} \varsigma}{1 + (\tau^2 + \tau^2\varsigma^2)^{\frac{(n-2s+2)2^*_{\sigma}}{2}}} \, d\varsigma  
= 
\frac{1}{2\tau^2} \int\limits_{\tau^2}^{+\infty} \frac{r^{ \frac{(1 - s + \sigma)2^*_{\sigma} - 2}{2}}}{1 + r^{\frac{(n-2s+2)2^*_{\sigma}}{2}}} \, dr 
\leq
\frac{1}{2\tau^2} \int\limits_{0}^{+\infty} \frac{r^{ \frac{(1 - s + \sigma)2^*_{\sigma} - 2}{2}}}{1 + r^{\frac{(n-2s+2)2^*_{\sigma}}{2}}} \, dr,
\end{equation*}
while for $\tau > 1$ we have
\begin{equation*}
\int\limits_{0}^{+\infty} \frac{(\tau^2 + \tau^2\varsigma^2)^{ \frac{(1 - s + \sigma)2^*_{\sigma} - 2}{2}} \varsigma}{1 + (\tau^2 + \tau^2\varsigma^2)^{\frac{(n-2s+2)2^*_{\sigma}}{2}}} \, d\varsigma  
\leq
\tau^{-(n-s-\sigma+1)2^*_{\sigma} - 2}
\int\limits_{0}^{+\infty} \frac{\varsigma}{(1+ \varsigma^2)^{\frac{(n-s-\sigma+1)2^*_{\sigma} + 2}{2}}} \, d\varsigma.
\end{equation*}
So, choosing sufficiently small $\beta,$ we get an estimate 
\begin{equation*}
P_{\varepsilon}(\tau) \leq C \left(\chi_{[0,1]}(\tau) \tau^{\alpha+n - 2 - \beta} + \chi_{[1, +\infty)}(\tau) \tau^{\alpha+n+\beta - (n-s-\sigma+1)2^*_{\sigma} - 2}\right)
\end{equation*}
with the summable majorant in the right-hand side (recall that $\alpha < n-2s+3$):
\begin{equation*}
\alpha+n+\beta - (n-s-\sigma+1)2^*_{\sigma} - 2 < -1 + \beta - \frac{2\sigma(n-2s+2)}{n - 2\sigma} < -1.
\end{equation*}
3. Using (\ref{f1_f_relations}), we obtain
\begin{multline*}
\left|\varepsilon \frac{I_4}{f(\varepsilon)} \right|
\leq 
O_{\delta}(1) \int\limits_{0}^{\frac{\delta}{\varepsilon}}  \frac{\tau^{n-2} f_1(\varepsilon \tau)}{\varepsilon |f(\varepsilon)|} \int\limits_{0}^{\sqrt{\frac{\delta^2}{\varepsilon^2} - \tau^2}} \frac{|\Phi(\tau, y_n)|^{2^*_{\sigma}} dy_n}{(\tau^2 + y_n^2)^{ \frac{(s-\sigma)2^*_{\sigma} + 2}{2}}} \, d\tau
\\ \leq 
o_{\delta}(1)  \int\limits_{0}^{\frac{\delta}{\varepsilon}}  \frac{\tau^{n} f(\varepsilon \tau)}{f(\varepsilon)}  \int\limits_{0}^{\sqrt{\frac{\delta^2}{\varepsilon^2 \tau^2} - 1}} \frac{|\Phi(\tau, \tau\varsigma)|^{2^*_{\sigma}} \, d\varsigma}{(\tau^2 + \tau^2\varsigma^2)^{ \frac{(s-\sigma)2^*_{\sigma} + 2}{2}}}  \, d\tau.
\end{multline*}
Similar to the previous estimate, the integral in the right-hand side has the finite limit as~$\varepsilon \to 0,$ what completes the proof. \qedhere
\end{proof}
To get (\ref{denom_est}) we put $\mathcal{A}_1 \left( \varepsilon \right) := I_3;$  estimates $I_4 = o_{\delta}(1) \mathcal{A}_1 \left( \varepsilon \right)$ and $I_2 = o_{\varepsilon}(1) \mathcal{A}_1 \left( \varepsilon \right)$ follow from Lemma \ref{lemma_Ik} and the inequality 
\begin{equation*}
I_2
\leq
C(\delta) \cdot \varepsilon^{\frac{n(n-2s+2)}{n-2\sigma}}
=
o_{\varepsilon}(1) \cdot \varepsilon^{\alpha - 1} 
\leq
o_{\varepsilon}(1) \cdot \varepsilon^{-1} f(\varepsilon)
=
o_{\varepsilon}(1) \cdot \mathcal{A}_1 \left( \varepsilon \right).
\end{equation*}
\section{Estimate of the numerator and derivation of (\ref{num_est})}
For brevity, we denote $\mathfrak{y} := \frac{n-2s}{2}.$ For $i \in [1:n-1]$ we have 
\begin{equation*}
\begin{pmatrix}
 \partial_{t} w_{\varepsilon}(X) \\
 \partial_{x_i} w_{\varepsilon}(X) \\
 \partial_{x_n} w_{\varepsilon}(X)
\end{pmatrix}
=
\begin{pmatrix}
\varepsilon^{-\mathfrak{y}-1} \mathcal{W}_{z}(\Theta_{\varepsilon}\left( X \right))\widetilde{\varphi}(x) \\[0.5em]
 \varepsilon^{-\mathfrak{y}-1} [\mathcal{W}_{y_i}(\Theta_{\varepsilon}(X)) - \mathcal{W}_{y_n}(\Theta_{\varepsilon}(X)) F_{x_i}(x')] \widetilde{\varphi}(x) +  \varepsilon^{-\mathfrak{y}} \mathcal{W}(\Theta_{\varepsilon}(X))\widetilde{\varphi}_{x_i}(x) \\[0.5em]
\varepsilon^{-\mathfrak{y}-1} \mathcal{W}_{y_n}(\Theta_{\varepsilon}(X))\widetilde{\varphi}(x) +  \varepsilon^{-\mathfrak{y}} \mathcal{W}(\Theta_{\varepsilon}(X))\widetilde{\varphi}_{x_n}(x)
\end{pmatrix}.
\end{equation*}
Using these formulae we get the representation for the energy 
\begin{multline*}
\mathcal{E}_s \left[ w_{\varepsilon} \right]
=
\int \limits_0^{+\infty} t^{1-2s} \int\limits_{\Omega} \Biggl( \sum_{i=1}^{n-1} \Bigl[\varepsilon^{-2\mathfrak{y} - 2} \widetilde{\varphi}^2 (x) \mathcal{W}^2_{y_i}(\Theta_{\varepsilon}(X))
\\-
2 \varepsilon^{-2\mathfrak{y}-2} \widetilde{\varphi}^2(x) \mathcal{W}_{y_i}(\Theta_{\varepsilon}(X)) \mathcal{W}_{y_n}(\Theta_{\varepsilon}(X)) F_{x_i} (x') 
+
2 \varepsilon^{-2\mathfrak{y}-1} \widetilde{\varphi}_{x_i}(x) \widetilde{\varphi} (x) \mathcal{W}_{y_i}(\Theta_{\varepsilon}(X)) \mathcal{W}(\Theta_{\varepsilon}(X))
\\ -
2\varepsilon^{-2\mathfrak{y}-1} \widetilde{\varphi}_{x_i}(x) \widetilde{\varphi}(x) F_{x_i}(x') \mathcal{W}_{y_n}(\Theta_{\varepsilon}(X))\mathcal{W}(\Theta_{\varepsilon}(X))
\\+
\varepsilon^{-2\mathfrak{y}-2} \widetilde{\varphi}^2(x) F^2_{x_i}(x')  \mathcal{W}_{y_n}^2(\Theta_{\varepsilon}(X))
+
\varepsilon^{-2\mathfrak{y}} \widetilde{\varphi}^2_{x_i}(x) \mathcal{W}^2(\Theta_{\varepsilon}(X)) \Bigr]
\\+
\varepsilon^{-2\mathfrak{y}-2} \widetilde{\varphi}^2(x) \mathcal{W}^2_{y_n}(\Theta_{\varepsilon}(X)) +
2\varepsilon^{-2\mathfrak{y}-1} \widetilde{\varphi}_{x_n}(x) \widetilde{\varphi}(x) \mathcal{W}_{y_n}(\Theta_{\varepsilon}(X)) \mathcal{W}(\Theta_{\varepsilon}(X)) 
\\+
\varepsilon^{-2\mathfrak{y}} \widetilde{\varphi}_{x_n}^2(x) \mathcal{W}^2(\Theta_{\varepsilon}(X))
+
 \varepsilon^{-2\mathfrak{y}-2} \widetilde{\varphi}^2(x) \mathcal{W}^2_{z}(\Theta_{\varepsilon}(X)) \Biggr) \, dX
 =:
 J_1-J_2+\dots+J_9+J_{10}
\end{multline*}
First, we estimate $J_1+ J_7 + J_{10}$ as follows:
\begin{equation*}
J_1+J_7 + J_{10} 
=
\int \limits_0^{+\infty}  z^{1-2s}  \int\limits_{\mathbb{R}^n_+} \varphi_{\delta}^2(\varepsilon y)|\nabla_{Y} \mathcal{W}(Y)|^2 \, dY 
= 
\mathcal{S}^{Sp}_{s, \sigma}(\mathbb{R}^n_+) -   \int\limits_{\mathbb{R}^n_+} \left[1-\varphi_{\delta}^2(\varepsilon y)\right] \cdot \mathcal{V}(y) \, dy.
\end{equation*}
From (\ref{GradW_common_ineq}) we get
\begin{equation*}
\int\limits_{\mathbb{R}^n_+} \left[1-\varphi_{\delta}^2(\varepsilon y)\right] \cdot \mathcal{V}(y) \, dy \leq C \int\limits_{\tfrac{\delta}{\varepsilon}}^{+\infty} r^{-3+2s-n} \, dr
= C \left(\frac{\varepsilon}{\delta}\right)^{n - 2s + 2}
\end{equation*}
what gives
\begin{equation*}
J_1+J_7 + J_{10} 
=
\mathcal{S}^{Sp}_{s, \sigma}(\mathbb{R}^n_+) + C(\delta) O(\varepsilon^{n - 2s + 2}).
\end{equation*}
Further, using (\ref{W_common_ineq}) and (\ref{GradW_common_ineq}) we estimate $J_3+J_8:$
\begin{multline*}
J_3+J_8 \leq 
2 \varepsilon \int \limits_0^{+\infty} z^{1-2s} \int\limits_{\mathbb{R}^n_+} \varphi_{\delta}(\varepsilon y) |\nabla_{y} \varphi_{\delta}(\varepsilon y)| \cdot \mathcal{W}(Y) |\nabla_{Y} \mathcal{W}(Y)| \, dY 
\\ \leq
\frac{C\varepsilon}{\delta} 
\left(\int\limits_{\mathbb{K}_{\frac{\delta}{2\varepsilon}}} \mathcal{V}(y) \, dy
\times
\int\limits_{\mathbb{K}_{\frac{\delta}{2\varepsilon}}} \int \limits_0^{+\infty} z^{1-2s} |\mathcal{W}(Y)|^2 \, dY \right)^{\frac{1}{2}}
\\ \leq
\frac{C\varepsilon}{\delta} \left( \ \int\limits_{\tfrac{\delta}{2\varepsilon}}^{\tfrac{\delta}{\varepsilon}} r^{-3+2s-n} \, dr 
\times
\int\limits_{\tfrac{\delta}{2\varepsilon}}^{\tfrac{\delta}{\varepsilon}} r^{-1+2s-n} \, dr \right)^{\frac{1}{2}}
= C \left(\frac{\varepsilon}{\delta}\right)^{n - 2s + 2}.
\end{multline*}
We estimate $J_4$ in a similar way:
\begin{multline*}
|J_4| \leq
2\varepsilon \int \limits_0^{+\infty}  z^{1-2s}  \int\limits_{\mathbb{R}^n_+} \varphi_{\delta}(\varepsilon y) |\nabla_{y} \varphi_{\delta}(\varepsilon y)| \mathcal{W}(Y) |\nabla_{Y} \mathcal{W}(Y)| |\nabla_{y'} F(\varepsilon y')| \, dY
\\ \leq
\frac{C\varepsilon}{\delta}  \int\limits_{\frac{\delta}{2\varepsilon}}^{\frac{\delta}{\varepsilon}} r^{2s-2n}  \int\limits_{0}^{r} \frac{\tau^{n-2} f_3(\varepsilon \tau) }{\sqrt{r^2 - \tau^2}} d\tau dr 
\leq
\frac{C\varepsilon^{n-2s+2}}{\delta}  \int\limits_{\frac{\delta}{2}}^{\delta} \tilde{r}^{2s-2n}  \int\limits_{0}^{\tilde{r}} \frac{\tilde{\tau}^{n-2} f_3(\tilde{\tau}) }{\sqrt{\tilde{r}^2 - \tilde{\tau}^2}} d\tilde{\tau} d\tilde{r} = C(\delta) \varepsilon^{n-2s+2}.
\end{multline*}
Also, (\ref{W_common_ineq}) allows to estimate $J_6+ J_9:$
\begin{equation*}
J_6+J_9
=
C \varepsilon^2 \int \limits_0^{+\infty} z^{1-2s} \int\limits_{\mathbb{R}^n_+} |\nabla_{y} \varphi_{\delta}(\varepsilon y)|^2 \mathcal{W}^2(Y) \, dY 
\leq 
 \frac{C \varepsilon^2}{\delta^2} \int\limits_{\tfrac{\delta}{2\varepsilon}}^{\tfrac{\delta}{\varepsilon}} r^{-1+2s-n} \, dr
=
C\left(\frac{\varepsilon}{\delta}\right)^{n - 2s + 2}.
\end{equation*}

Now we transform the main term $J_2.$ Integrating by parts, we obtain
\begin{multline*}
J_2 = 
-\frac{2}{\varepsilon} \int \limits_0^{+\infty} z^{1-2s} \int\limits_{\mathbb{R}^n_+} \sum_{i=1}^{n-1} \Bigl[ \varphi_{\delta}^2(\varepsilon y) \mathcal{W}_{y_i y_i}(Y) \mathcal{W}_{y_n}(Y)F(\varepsilon y')
+
\left[\varphi_{\delta}^2(\varepsilon y)\right]_{y_i} \mathcal{W}_{y_i}(Y) \mathcal{W}_{y_n}(Y)F(\varepsilon y')
\\+ 
\varphi_{\delta}^2(\varepsilon y) \mathcal{W}_{y_i}(Y) \mathcal{W}_{y_iy_n}(Y)F(\varepsilon y') \Bigr] \, dY.
\end{multline*}
Next, we use the BVP (\ref{eq_st_tor}) to express the sum of second derivatives:
\begin{multline*}
J_2 = \frac{2}{\varepsilon} \int \limits_0^{+\infty} \int\limits_{\mathbb{R}^n_+}  \varphi_{\delta}^2(\varepsilon y) \mathcal{W}_{y_n}(Y)\left[ z^{1-2s}\mathcal{W}_{y_n y_n}(Y) + [z^{1-2s}\mathcal{W}_{z}(Y)]_{z} \right] F(\varepsilon y') \, dY 
\\-
\frac{2}{\varepsilon} \int \limits_0^{+\infty} \int\limits_{\mathbb{R}^n_+} z^{1-2s} \sum_{i=1}^{n-1}  \left[\varphi_{\delta}^2(\varepsilon y)\right]_{y_i} \mathcal{W}_{y_i}(Y)\mathcal{W}_{y_n}(Y)F(\varepsilon y')  \, dY
\\-
\frac{1}{\varepsilon} \int \limits_0^{+\infty} \int\limits_{\mathbb{R}^n_+}  z^{1-2s} \varphi_{\delta}^2(\varepsilon y) \left[ |\nabla_{y'}\mathcal{W}(Y)|^2 \right]_{y_n} F(\varepsilon y') \, dY
 =: \mathcal{H} + E_1 + E_2.
\end{multline*}
Intergrating by parts once more, we transform $\mathcal{H}$ as follows:
\begin{multline*}
\mathcal{H} 
=
\frac{1}{\varepsilon} \int \limits_0^{+\infty} \int\limits_{\mathbb{R}^n_+}  \varphi_{\delta}^2(\varepsilon y) F(\varepsilon y') \left[ z^{1-2s} \left[\mathcal{W}^2_{y_n}(Y)\right]_{y_n} + 2  \mathcal{W}_{y_n}(Y) \left[z^{1-2s}\mathcal{W}_{z}(Y)\right]_{z}\right] dY
\\= 
-\frac{1}{\varepsilon} \int \limits_0^{+\infty} z^{1-2s}  \Bigl[ \int\limits_{\mathbb{R}^{n-1}} \varphi^2_{\delta}(\varepsilon y') \mathcal{W}^2_{y_n}(y', 0, z) F(\varepsilon y') dy' + \int\limits_{\mathbb{R}^n_+} [\varphi_{\delta}^2(\varepsilon y)]_{y_n} \mathcal{W}^2_{y_n}(Y) F(\varepsilon y') \, dy \Bigr] dz 
\\+
\frac{2\mathcal{S}^{Sp}_{s, \sigma}(\mathbb{R}^n_+)}{\varepsilon} \int\limits_{\mathbb{R}^n_+}  \varphi_{\delta}^2(\varepsilon y) F(\varepsilon y') \Phi_{y_n}(y)  \frac{\Phi^{2^*_{\sigma} - 1}(y)}{|y|^{(s-\sigma)2^*_{\sigma}}}  \, dy
-
\frac{1}{\varepsilon} \int \limits_0^{+\infty} \int\limits_{\mathbb{R}^n_+}  \varphi_{\delta}^2(\varepsilon y) F(\varepsilon y') \left[z^{1-2s}\mathcal{W}^2_{z}(Y)\right]_{y_n} dY
\\=:
-E_3 + E_4 + \mathcal{K} + E_7.
\end{multline*}
We integrate by parts $\mathcal{K}$ and $E_7,$ taking into account $\mathcal{W}_{z}(y', 0, z) = 0,$ and obtain
\begin{gather*}
\mathcal{K} = \frac{2\mathcal{S}^{Sp}_{s, \sigma}(\mathbb{R}^n_+)}{\varepsilon} \int\limits_{\mathbb{R}^n_+}  \left[ -\frac{\left[\varphi_{\delta}^2(\varepsilon y) \right]_{y_n}}{2^*_{\sigma}} + \varphi_{\delta}^2(\varepsilon y) \frac{(s-\sigma) y_n}{|y|^2}  \right]  \frac{\Phi^{2^*_{\sigma}}(y)}{|y|^{(s-\sigma)2^*_{\sigma}}} F(\varepsilon y') \, dy 
=: E_5 + E_6, \\
E_7 = \frac{1}{\varepsilon} \int \limits_0^{+\infty} z^{1-2s} \int\limits_{\mathbb{R}^n_+}  \left[\varphi_{\delta}^2(\varepsilon y) \right]_{y_n}F(\varepsilon y') \mathcal{W}^2_{z}(Y) \, dY.
\end{gather*}
\begin{lemma}
\label{lemma_Jk}
The following relations hold:
\begin{enumerate}
\item $|E_1 + E_2 + E_4 + E_7| = C(\delta) \varepsilon^{n -2s + 2};$
\item \begin{equation}
\label{E3_estim}
\lim\limits_{\varepsilon \to 0} \varepsilon \frac{E_3}{f(\varepsilon)}  = C \int \limits_{0}^{+\infty} \tau^{n+\alpha-2} \int \limits_0^{+\infty} z^{1-2s} |\nabla_{\tau,z} \mathcal{W}(\tau, 0, z)|^2 \, dzd\tau < +\infty;
\end{equation}
\item $|E_5|= o(\varepsilon^{n -2s + 2});$
\item $E_6 =  \frac{2\mathcal{S}^{Sp}_{s, \sigma}(\mathbb{R}^n_+)}{2^*_{\sigma}} \mathcal{A}_1 \left( \varepsilon \right) \cdot (1 + o_{\varepsilon}(1)).$
\end{enumerate}
\end{lemma}
\begin{proof}
1. The statement follows from the following inequalities
\begin{multline*}
|E_1 + E_2 + E_4 + E_7| 
= 
\left| \frac{1}{\varepsilon} \int \limits_0^{+\infty} z^{1-2s} \int\limits_0^{+\infty} \int\limits_0^{+\infty} \Bigl[ -2[\varphi_{\delta}^2(\varepsilon y)]_{\tau} \mathcal{W}_{\tau}(Y)\mathcal{W}_{y_n}(Y) \right.
\\+ 
\left. [\varphi_{\delta}^2(\varepsilon y)]_{y_n} \left[\mathcal{W}^2_{\tau}(Y) - \mathcal{W}^2_{y_n}(Y) + \mathcal{W}^2_{z}(Y) \right] \Bigr] \, dy_n
\times \int\limits_{\mathbb{S}^{n-2}_{\tau}} F(\varepsilon y') \, d\mathbb{S}^{n-2}_{\tau}(y') \, d\tau dz\right| 
\\ \leq
\frac{C}{\delta} \int\limits_{\frac{\delta}{2\varepsilon}}^{\frac{\delta}{\varepsilon}} \int\limits_0^{r} \int \limits_0^{+\infty} z^{1-2s} \left|\nabla_{\tau, y_n, z} \mathcal{W}(\tau, \sqrt{r^2 - \tau^2}, z)\right|^2  dz \frac{r\tau^{n-2}}{\sqrt{r^2 - \tau^2}} |f(\varepsilon \tau)| \, d\tau dr
\\ \leq
\frac{C}{\delta} \int\limits_{\frac{\delta}{2\varepsilon}}^{\frac{\delta}{\varepsilon}} r^{-2n + 2s - 1} \int\limits_0^{r} \frac{\tau^{n-2}|f(\varepsilon \tau)|}{\sqrt{r^2 - \tau^2}} \, d\tau dr 
= 
C(\delta) \varepsilon^{n-2s+2} \int\limits_{\frac{\delta}{2}}^{\delta} \int\limits_0^{\tilde{r}} \frac{\tilde{r}^{-2n + 2s - 1} \tilde{\tau}^{n-2}|f(\tilde{\tau})|}{\sqrt{\tilde{r}^2 - \tilde{\tau}^2}} \, d\tilde{\tau} d\tilde{r}.
\end{multline*}

2. As for the estimate of $I_3$ in Lemma \ref{lemma_Ik}, we use the Lebesgue theorem: since
\begin{equation*}
\varepsilon \frac{E_3}{f(\varepsilon)} = C \int\limits_{0}^{+\infty} \tau^{n-2} \frac{f(\varepsilon \tau)}{f(\varepsilon)} \varphi_{\delta}^2(\varepsilon \tau) \int \limits_0^{+\infty} z^{1-2s} |\nabla_{\tau,z} \mathcal{W}(\tau, 0, z)|^2 \, dzd\tau=:
C\int\limits_{0}^{+\infty} Q_{\varepsilon}(\tau) \, d\tau,
\end{equation*}
we get the integrand in the right-hand side of (\ref{E3_estim}) as the pointwise limit of $Q_{\varepsilon}(\tau).$ To construct the majorant we use (\ref{GradW_common_ineq}) and (\ref{psi_equal}):
\begin{multline*}
Q_{\varepsilon}(\tau)
\leq 
 \chi_{\left[0, \frac{\delta}{\varepsilon}\right]}(\tau) \tau^{n-2 + \alpha} \frac{\psi(\varepsilon \tau)}{\psi(\varepsilon)}  \int \limits_0^{+\infty} z^{1-2s} |\nabla_{\tau,z} \mathcal{W}(\tau, 0, z)|^2 \, dz
\leq 
\chi_{\left[0, \frac{\delta}{\varepsilon}\right]}(\tau)  \frac{\psi(\varepsilon \tau)}{\psi(\varepsilon)}  \frac{C \tau^{n-2 + \alpha}}{1 + |\tau|^{2n-2s+2}}
\\  \leq
C(\delta) \left(\chi_{[0,1]}(\tau) \cdot \tau^{\alpha+n - \beta - 2} + \chi_{[1, +\infty)}(\tau) \cdot \tau^{\alpha- n - 4 + \beta + 2s}\right),
\end{multline*}
what is summable for sufficiently small $\beta$ due to $\alpha < n - 2s + 3.$

3. We have:
\begin{multline*}
|E_5| 
=
\left|\frac{2\mathcal{S}^{Sp}_{s, \sigma}(\mathbb{R}^n_+)}{ 2^*_{\sigma} \cdot \varepsilon } \int\limits_0^{+\infty}  \int\limits_{\mathbb{S}^{n-2}_{\tau}} F(\varepsilon y') \, d\mathbb{S}^{n-2}_{\tau}(y')  \int\limits_0^{+\infty} \left[\varphi_{\delta}^2(\varepsilon \sqrt{\tau^2 + y_n^2}) \right]_{y_n} \frac{|\Phi|^{2^*_{\sigma}}(\tau, y_n)}{|y|^{(s-\sigma)2^*_{\sigma}}} \, dy_n d\tau \right|
\\ \leq 
\frac{C}{\delta} \int\limits_{\frac{\delta}{2\varepsilon}}^{\frac{\delta}{\varepsilon}} \int\limits_0^{r} \frac{|\Phi|^{2^*_{\sigma}}(\tau, \sqrt{r^2 - \tau^2})}{r^{(s-\sigma)2^*_{\sigma}}} \frac{r \tau^{n-2}|f(\varepsilon \tau)|}{\sqrt{r^2 - \tau^2}} \, d\tau dr
\leq
\frac{C}{\delta} \int\limits_{\frac{\delta}{2\varepsilon}}^{\frac{\delta}{\varepsilon}} \int\limits_0^{r} \frac{r\tau^{n-2}|f(\varepsilon \tau)|}{r^{2^*_{\sigma}(n -s-\sigma+1)} \sqrt{r^2 - \tau^2}} \, d\tau dr
\\ \leq
\frac{C \varepsilon^{2^*_{\sigma}(n -s-\sigma+1) - n}}{\delta} \int\limits_{\frac{\delta}{2}}^{\delta} \tilde{r}^{1 - 2^*_{\sigma}(n -s-\sigma+1)} \int\limits_0^{\tilde{r}} \frac{\tilde{\tau}^{n-2}|f(\tilde{\tau})|}{\sqrt{\tilde{r}^2 - \tilde{\tau}^2}} \, d\tau d\tilde{r} = o(\varepsilon^{n-2s+2}).
 \end{multline*}

4. Notice that the expression for $E_6$ coincides with the expression for $I_3$ up to two differences: 
we replace $\varphi_{\delta}^{2^*_{\sigma}}(\varepsilon y)$ with $\varphi_{\delta}^2(\varepsilon y)$ and multiply by $ \frac{2\mathcal{S}^{Sp}_{s, \sigma}(\mathbb{R}^n_+)}{2^*_{\sigma}}.$ Thus the statement follows from the argument from Lemma \ref{lemma_Ik}. \qedhere
\end{proof}
Lemma \ref{lemma_Jk} together with estimates $I_3 \asymp f(\varepsilon) \varepsilon^{-1} \asymp E_3$ and $\varepsilon^{n-2s+2} = o(f(\varepsilon) \varepsilon^{-1})$ gives
\begin{equation*}
J_2 = -E_3 \cdot (1 + o_{\delta}(1) + o_{\varepsilon}(1)) + \frac{2\mathcal{S}^{Sp}_{s, \sigma}(\mathbb{R}^n_+)}{2^*_{\sigma}}  \mathcal{A}_1 \left( \varepsilon \right) \cdot (1 + o_{\varepsilon}(1)).
\end{equation*}
It remains to estimate $J_5.$ Using (\ref{GradW_common_ineq}) and (\ref{f1_f_relations}) we get
\begin{multline*}
J_5 =   
\int \limits_0^{+\infty}  z^{1-2s}  \int\limits_{\mathbb{R}^n_+} \varphi_{\delta}^2(\varepsilon y) |\nabla_{y'} F(\varepsilon y')|^2 \mathcal{W}^2_{y_n}(Y) \, dY
\leq
C  \int\limits_{0}^{\tfrac{\delta}{\varepsilon}} \tau^{n-2} f_2(\varepsilon \tau) \, d\tau \int\limits_{0}^{\sqrt{\tfrac{\delta^2}{\varepsilon^2} - \tau^2}}  \mathcal{V}(\tau, y_n) \, dy_n 
\\ \leq 
\int\limits_{0}^{\tfrac{\delta}{\varepsilon}} \int\limits_{0}^{+\infty} \frac{C \tau^{n-2} f_2(\varepsilon \tau) }{(1 + \tau^2 + y^2_n)^{n-s+1}} \, dy_nd\tau
\leq
\int\limits_{0}^{\tfrac{\delta}{\varepsilon}} \frac{C \tau^{n-2} f_2(\varepsilon \tau)}{(1+\tau^2)^{\frac{2n-2s+ 1}{2}}} \, d\tau = \frac{o_{\delta}(1)}{\varepsilon} \int\limits_{0}^{\tfrac{\delta}{\varepsilon}} \frac{\tau^{n-3} |f(\varepsilon \tau)|}{(1+\tau^2)^{\frac{2n-2s+ 1}{2}}} \, d\tau.
\end{multline*}
The last integral can be estimated in the same way as $E_3$ in Lemma \ref{lemma_Jk}. This estimate gives~$J_5 = o_{\delta}(1) E_3.$ 

Denoting $\mathcal{A}_2 (\varepsilon) := E_3,$ we obtain (\ref{num_est}).


\begin{thebibliography}{CS}
\small
\bibitem{Aubin}
T. Aubin, {\it Problemes isop{\'e}rim{\'e}triques et espaces de Sobolev}, J. Diff. Geom., {\bf 11} (1976), no.~4, 573-598.

\bibitem{Caffarelli}
L. Caffarelli\ and\ L. Silvestre, {\it An extension problem related to the fractional Laplacian}, Comm. Part. Diff. Eqs., {\bf 32} (2007), no.~7-9, 1245-1260.

\bibitem{Cabre}
X. Cabr{\'e}\ and\ Y. Sire, {\it Nonlinear equations for fractional Laplacians, I: Regularity, maximum principles and Hamiltonian estimates}, Ann. Inst. H. Poincar{\'e}. Anal. Nonlin., {\bf 31} (2014), no.~1, 23-53.

\bibitem{Capella}
A. Capella, J. D{\'a}vila, L. Dupaigne, and Y. Sire, {\it Regularity of radial extremal solutions for some non-local semilinear equations}, Comm. Part. Diff. Eqs., {\bf 36} (2011), no.~8, 1353-1384.

\bibitem{Tav}
A. Cotsiolis and N.K. Tavoularis, {\it Best constants for Sobolev inequalities for higher order fractional derivatives}, J. Math. Anal. Appl., {\bf 295} (2004), no.~1, 225-236.

\bibitem{Dem_Naz}
A.V. Demyanov and A.I. Nazarov, {\it On solvability of the Dirichlet problem to the semilinear Schr{\"o}dinger equation with singular potential},  ZNS  POMI, {\bf 336} (2006), 25–45, (Russian); English transl.: J. Math. Sci., {\bf 143} (2007), no.~2, 2857-2868.

\bibitem{Egn}
H. Egnell, {\it Positive solutions of semilinear equations in cones}, Trans. Amer. Math. Soc., {\bf 330}~(1992), no.~1, 191-201.

\bibitem{Fall}
M.M. Fall and T. Weth, {\it Nonexistence results for a class of fractional elliptic boundary value problems}, J. Func. Anal., {\bf 263}~(2012), no.~8, 2205-2227.

\bibitem{Gh_Kang}
N. Ghoussoub and X.S. Kang, {\it Hardy–Sobolev critical elliptic equations with boundary singularities}, Ann. Inst. H. Poincar{\'e}. Anal. Nonlin., {\bf 21} (2004), no.~6, 767-793.

\bibitem{Gh_Rob}
N. Ghoussoub and F. Robert, {\it The effect of curvature on the best constant in the Hardy–Sobolev inequalities}, GAFA, {\bf 16} (2006), no.~6, 1201-1245.

\bibitem{Gh_Yu}
N. Ghoussoub and C. Yuan, {\it Multiple solutions for quasilinear PDEs involving the critical Sobolev and Hardy exponents}, Trans. Amer. Math. Soc., {\bf 352} (2000), no.~12, 5703-5743.

\bibitem{Glaser}
V. Glaser, W.E. Thirring, Н. Grosse, and A. Martin, {\it A family of optimal conditions for the absence of bound states in a potential}, Les rencontres physiciens-math{\'e}maticiens de Strasbourg, {\bf 23} (1976), no.~1, 0-21.

\bibitem{Hardy}
G. Hardy, J.E. Littlewood and G. Polya, {\it Inequalities}, Cambridge Univ. Press, 1934.

\bibitem{Herbst}
I.W. Herbst, {\it Spectral theory of the operator $(p^2 + m^2)^{1/2} - Z e^2/r$}, Comm. Math. Phys., {\bf 53}~(1977), no.~3, 285-294.

\bibitem{Ilin}
V.P. Il'in, {\it Some integral inequalities and their applications in the theory of differentiable functions of several variables}, Mat. Sb., {\bf 54} (1961), no.~3, 331-380 (Russian).
 
\bibitem{Kawohl}
B. Kawohl, {\it Rearrangements and convexity of level sets in PDE}, Springer Lecture Notes in Math., {\bf 1150} (1985).

\bibitem{Lad_Ur}
O.A. Ladyzhenskaya and N.N. Ural’tseva, {\it Linear and Quasilinear Equations of Elliptic Type}, 2nd ed., Nauka, 1973, (Russian).

\bibitem{Lieb}
E. Lieb, {\it Sharp constants in the Hardy-Littlewood-Sobolev and related inequali­ties}, Ann. Math., \textbf{118} (1983), no.~2, 349-374.

\bibitem{Lieb_Lo}
E. Lieb and M. Loss, {\it Analysis}, Grad. Studies Math., 2nd ed., Amer. Math. Soc. {\bf 14} (2001).

\bibitem{Lions}
P.L. Lions, {\it The concentration-compactness principle in the Calculus of Variations. The locally compact case},  Ann. Inst. H. Poincar{\'e}. Anal. Nonlin., {\bf 1} (1984), 109-145, 223-283. {\it The limit case}, Rev. Mat. Iberoam. {\bf 1} (1985), 45-121, 145-201.

\bibitem{MN1}
R. Musina and A.I. Nazarov, {\it Fractional Hardy-Sobolev inequalities on half spaces}, Nonlin. Analysis – TMA, {\bf 178} (2019), 32-40.

\bibitem{MN2}
R. Musina and A.I. Nazarov, {\it On fractional Laplacians}, Comm. Part. Diff. Eqs., {\bf 39} (2014), no.~9, 1780-1790.

\bibitem{MN3}
R. Musina and A.I. Nazarov, {\it On fractional Laplacians--3}, ESAIM: COCV, {\bf 22} (2016), no.~3, 832-841.

\bibitem{MN4}
R. Musina and A.I. Nazarov, {\it On the Sobolev and Hardy constants for the fractional Navier Laplacian}, Nonlin. Analysis -- TMA, {\bf 121} (2015), 123-129.

\bibitem{MN5}
R. Musina and A.I. Nazarov, {\it Sobolev inequalities for fractional Laplacians on half spaces}, Adv. Calc. Var. (2018), DOI: https://doi.org/10.1515/acv-2018-0020.

\bibitem{MN6}
R. Musina and A.I. Nazarov, {\it Strong maximum principles for fractional Laplacians}, Proc. Roy. Soc. Edinburgh Sect. A, (2019), 1-18, DOI: https://doi.org/10.1017/prm.2018.81.

\bibitem{Naz_Cone}
A.I. Nazarov, {\it Hardy–Sobolev inequalities in a cone}, Probl. Math. Anal., {\bf 31} (2005), 39-46 (Russian); English transl.: J. Math. Sci. {\bf 132} (2006), no.~4, 419-427.

\bibitem{Ros_Oton}
X. Ros-Oton and J. Serra, {\it The Pohozaev identity for the fractional Laplacian}, ARMA, {\bf 213} (2014), no.~2, 587-628.

\bibitem{Seneta}
E. Seneta, {\it Regularly varying functions}, Lect. Notes Math., Springer Verlag, {\bf 508} (1976).

\bibitem{Scheglova}
A.P. Shcheglova, {\it The Neumann boundary value problem for a semilinear elliptic equation in a thin cylinder. The least energy solutions},  ZNS  POMI, {\bf 348} (2007), 272–302, (Russian); English transl.: J. Math. Sci., {\bf 152} (2008), no.~5, 780-798.

\bibitem{Stinga}
P.R. Stinga and J.L. Torrea, {\it Extension problem and Harnack's inequality for some fractional operators}, Comm. Part. Diff. Eqs., {\bf 35} (2010), no.~11, 2092-2122.

\bibitem{Talenti}
G. Talenti, {\it Best constant in Sobolev inequality}, Ann. di Mat. Pura ed Appl., {\bf 110} (1976), no.~1, 353-372.

\bibitem{Triebel}
H. Triebel, {\it Interpolation theory, function spaces, differential operators}, Deutscher Verlag Wissensch., Berlin, 1978.

\bibitem{Ustinov}
N.S. Ustinov, {\it On attainability of the best constant in fractional Hardy--Sobolev inequality with the Spectral Dirichlet Laplacian}, Funct. An. and Appl., {\bf 53} (2019), no.~3 (to appear).

\bibitem{Yang}
J. Yang, {\it Fractional Sobolev--Hardy inequality in $\mathbb{R}^n$},  Nonlin. Analysis -- TMA, {\bf 119} (2015), 179-185.

\end{thebibliography}
\end{document}